\documentclass{article}

\usepackage[english]{babel}
\usepackage[table,xcdraw]{xcolor}
\usepackage{algorithm}
\usepackage{algpseudocode}
\usepackage{forest}
\usepackage{comment}
\usepackage{lscape} 
\usepackage{amsfonts}
\usepackage[style=ieee]{biblatex}
\usepackage{optidef}
\usepackage{tikz}
\usepackage[verbose]{placeins} 
\usepackage{nicematrix}
\usepackage{optidef}
\usepackage{amsthm}
\usetikzlibrary{decorations.pathreplacing} % tikz library should be here, in preamble

\pgfdeclarelayer{background}
\pgfdeclarelayer{foreground}
\pgfsetlayers{background,main,foreground}

\tikzstyle{standard}=[fill=white, draw=black, shape=circle]
\tikzstyle{new style 0}=[fill={rgb,255: red,255; green,128; blue,0}, draw=black, shape=circle]
\tikzstyle{new style 1}=[fill=blue, draw=black, shape=circle]

\tikzstyle{new edge style 0}=[-, draw={rgb,255: red,191; green,191; blue,191}, line width=1pt]
\tikzstyle{new edge style 1}=[->]
\tikzstyle{new edge style 2}=[-, draw={rgb,255: red,255; green,128; blue,0}]
\tikzstyle{new edge style 3}=[-]

\tikzset{
  tree node/.style = {align=center, inner sep=0pt, font = \scriptsize},
  S/.style = {draw, circle, minimum size = 8mm, top color=white, bottom color=white!20},
  tree node label/.style={font=\scriptsize},
}
\forestset{
  declare toks={left branch prefix}{},
  declare toks={right branch prefix}{},
  declare toks={left branch suffix}{},
  declare toks={right branch suffix}{},
  tree node left label/.style={
    label=170:#1,
  },
  tree node right label/.style={
    label=10:#1,
  },
  maths branch labels/.style={
    branch label/.style={
      if n=1{
        edge label={node [left, midway] {$\forestoption{left branch prefix}##1\forestoption{left branch suffix}$}},
      }{
        edge label={node [right, midway] {$\forestoption{right branch prefix}##1\forestoption{right branch suffix}$}},
      }
    },
  },
  text branch labels/.style={
    branch label/.style={
      if n=1{
        edge label={node [left, midway] {\foresteoption{left branch prefix}##1\forestoption{left branch suffix}}},
      }{
        edge label={node [right, midway] {\forestoption{right branch prefix}##1\forestoption{right branch suffix}}},
      }
    },
  },
  text branch labels,
  set branch labels/.style n args=4{%
    left branch prefix={#1},
    left branch suffix={#2},
    right branch prefix={#3},
    right branch suffix={#4},
  },
  set maths branch labels/.style n args=4{
    maths branch labels,
    set branch labels={#1}{#2}{#3}{#4},
  },
  set text branch labels/.style n args=4{
    text branch labels,
    set branch labels={#1}{#2}{#3}{#4},
  },
  branch and bound/.style={
    /tikz/every label/.append style=tree node label,
    maths branch labels,
    for tree={
      tree node,
      S,
      math content,
      s sep'+=20mm,
      l sep'+=5mm,
      thick,
      edge+={thick},
    },
    before typesetting nodes={
      for tree={
        split option={content}{:}{tree node left label,content,tree node right label,branch label},
      },
    },
    where n children=0{
    }{},
  },
  }
\usepackage[letterpaper,top=2cm,bottom=2cm,left=3cm,right=3cm,marginparwidth=1.75cm]{geometry}

\usepackage{amsmath}
\usepackage{graphicx}
\usepackage[colorlinks=true, allcolors=blue]{hyperref}
\newtheorem{theorem}{Theorem}
\numberwithin{theorem}{section}
\newtheorem{definition}{Definition}
\numberwithin{definition}{section}
\newtheorem{lemma}{Lemma}
\numberwithin{lemma}{section}
\newtheorem{conjecture}{Conjecture}
\numberwithin{conjecture}{section}
\newtheorem{proposition}{Proposition}
\numberwithin{proposition}{section}

\title{Constructing extremal triangle-free graphs using integer programming}
\author{Ali Erdem Banak, Tınaz Ekim, Z. Caner Taşkın\\ \normalsize{Bo\u{g}aziçi University, Department of Industrial Engineering}}

\begin{document}
\emergencystretch 3em
\maketitle

\begin{abstract}
The maximum number of edges in a graph with matching number $m$ and maximum degree $d$ has been determined in \cite{ch76} and \cite{bk09}, where some extremal graphs have also been provided. Then, a new question has emerged: how the maximum edge count is affected by forbidding some subgraphs occurring in these extremal graphs? In \cite{aey22}, the problem is solved in triangle-free graphs for $d \ge m$, and for $d < m$ with either $Z(d) \le m < 2d$ or $d \leq 6$, where $Z(d)$ is approximately $5d/4$. The authors derived structural properties of triangle-free extremal graphs, which allows us to focus on constructing small extremal components to form an extremal graph. Based on these findings, in this paper, we develop an integer programming formulation for constructing extremal graphs. Since our formulation is highly symmetric, we use our own implementation of Orbital Branching to reduce symmetry. We also implement our integer programming formulation so that the feasible region is restricted iteratively. Using a combination of the two approaches, we expand the solution into $d \leq 10$ instead of $d \leq 6$ for $m > d$. Our results endorse the formula for the number of edges in all extremal triangle-free graphs conjectured in \cite{aey22}.\\

\end{abstract}

\section{Introduction}

Extremal graph theory studies how big or small a graph parameter (usually the number of edges or vertices) can be under some local constraints.  A graph that is an optimal solution of such a maximization or minimization problem is called an \textit{extremal graph}. Turan's graphs are one of the best known examples of extremal graph theory results \cite{e64}. The extremal problem of maximizing the number of edges in a graph while limiting the size of a maximum matching and the maximum degree has been first posed by Erd\"os and Rado \cite{er60} in a more general context. This question has been first answered by Chvátal and Hanson in \cite{ch76} using some optimization techniques and Berge's matching formula. Later, Balachandran and Khare \cite{bk09} provided a constructive proof for the same problem. The authors describe some extremal graphs having three types of components, namely star graphs, complete graphs, and almost complete graphs (which contain $C_4$'s, that is induced cycles of length 4). This construction paves the way for new variants of the problem; how the maximum number of edges will be affected if we exclude one type of component described in \cite{bk09}? Accordingly, one can consider the same extremal problem for graphs (independently) restricted to be claw-free, which forbids the smallest star graph; triangle-free, which forbids the smallest complete graph; or $C_4$-free, which forbids almost-complete graphs.

Dibek, Ekim, and Heggernes studied claw-free graphs and settled the discussion by showing the cases where the maximum number of edges achieves the maximum number for general graphs and where it is strictly less (than for the general case) \cite{deh17}. Blair, Heggernes, Lima, and Lokshtanov worked on chordal graphs, which exclude $C_4$'s and therefore almost-complete graphs. They show that the maximum edge count in chordal graphs is the same as the general graphs \cite{b22}. Similarly,  Maland \cite{m15} answers the question on bipartite graphs, split graphs, disjoint union of split graphs, and unit interval graphs.

Triangle-free graphs have been considered from the same perspective by Ahanjideh, Ekim, and Y\i ld\i z in \cite{aey22}. In their study, the authors provided some partial answers to the extremal problem, which is thus not completely solved yet. The authors constructed extremal graphs for all $d \ge m$. For the cases with $d < m$, they solved the problem for either $d \le 6$ or $Z(d) \le m < 2d$ where $Z(d)$ is approximately $5d/4$. Consequently, the cases with $d \ge 7$ for either $m\geq 2d$ or $d < m < Z(d)$ are left open in \cite{aey22}. Nevertheless, the authors derived some structural results that shed light on these open cases: For all natural numbers $d\geq 2$ and $m$, there exists an extremal graph whose components are either star graphs or factor-critical triangle-free extremal graphs with maximum degree $d$ and matching number between $d$ and $Z(d)$. This result implies that to construct extremal graphs for $2d \le m$, it is sufficient to find extremal components, that is, extremal graphs with matching number between $d$ and $Z(d)$ and decide how many of each one of those components and how many star graphs should be taken. This problem is then expressed as a Knapsack Formulation where the total volume is $m$, individual volumes are the matching numbers of each extremal component, and utilities are maximum edge counts in each extremal component. 

Following the above observation, in this paper, we focus on the computation of maximum edge counts in extremal components and propose integer programming approaches to find them. To the best of our knowledge, the use of integer programming formulations in the construction of extremal graphs is rather new. A few exceptions are the recent paper on Ramsey numbers to improve the known lower bounds \cite{furini} and the construction of (small) extremal chemical graphs with given number of vertices and/or degrees optimizing some invariants \cite{chem1,chem2}. In Section 2, we introduce the notation and mention the known results for general graphs and triangle-free graphs. We also provide the knapsack formulation in \cite{aey22} to compute the max edge counts for $d \ge 7$ and $m \ge 2d$. We state Conjecture  \ref{conj:comp} and Conjecture \ref{conj:lemma}, giving the formula for the max edge counts of extremal components and all triangle-free extremal graphs respectively. In particular, we exhibit the results in \cite{aey22}, which allows us to restrict our study to the construction of extremal factor-critical triangle-free components with $d > 6$ and $d < m < Z(d)$. In Section 3, we explain our methodology to construct the extremal components via four different exact methods. Our Basic Formulation exploits the structural information provided in \cite{aey22} to bound the matching number without explicitly including it in the formulation. Noting that the Basic Formulation is highly symmetric, we propose various methods to reduce symmetry. In particular, we utilize Orbital Branching \cite{o11}, which identifies equivalent variables and sets the value of more than one of them in a branch to reduce the symmetry within the branch-and-cut algorithm. Next, Iterative Method solves the Basic Formulation by considering only a portion of the feasible region at each iteration. In this approach, we set all the degrees of vertices to $d$ and check if the problem is feasible. If not, we iteratively check the existence of graphs with the next possible upper bound for the edge count. Finally, we combine both methods, applying Orbital Branching within the Iterative Method. In Section 4, we discuss our implementation, where we used CPLEX 20.1.0 for integer programming and nauty for calculating the orbits. 

In Section 5, we provide our computational results. We show that the Orbital Branching and the Iterative Method are both better than the Basic Formulation using the default branching strategy; moreover, their combination gives us the best results. Our method allows us to construct all extremal components for $d=7,8,9,10$ and $d< m < Z(d)$. Subsequently, we use the maximum edge counts for these extremal components as the parameters of the Knapsack Formulation proposed in \cite{aey22}. We solve the Knapsack Formulation for $d=7,8,9,10$ within $0.2$ seconds even for $m = 2000$. Therefore, the edge counts of edge-extremal graphs for all possible $m$ values are found for $d=7,8,9,10$. We could also construct some (but not all) extremal components for $d=11,12$ and 13 with $d<m<Z(d)$. Our findings support the conjectures suggested by Ahanjideh, Ekim, and Y\i ld\i z in \cite{aey22} and explained in Section \ref{sec:prem}. In particular, we now have a stronger evidence that the formula provided in \cite{aey22} (see Theorem \ref{thm:main} and Conjecture \ref{conj:lemma}) gives the maximum number of edges in a triangle-free graph with maximum degree at most $d$ and matching number at most $m$ for all natural numbers $d \ge 2$ and $m$. Yet, a formal proof remains as a future research.

\section{Notation and Preliminaries}\label{sec:prem}

In this paper, undirected graphs are represented with ${G = (V(G), E(G))}$ where ${V(G)}$ is the set of vertices and ${E(G)}$ is the set of edges of the graph $G$. The {\em degree} of a vertex $v$, denoted by ${d(v)}$, is the number of vertices adjacent to $v$. The {\em maximum degree} of a graph $G$ is the maximum degree of a vertex in $G$, denoted by $\Delta(G)$. A {\em matching} of a graph $G$ is defined as a set of edges with no common vertices. The size of a maximum matching of a graph $G$ is denoted by ${\nu(G)}$. 

Maximum edge count for a general graph with maximum degree at most $d$ and matching number at most $m$ is shown with $f_{GEN}(d, m)$. It can be observed that if we do not limit either $\Delta(G)$ or $\nu(G)$, a graph can have an unlimited number of edges; such examples are a central vertex with an unbounded number of neighbors, and an unbounded number of independent edges respectively. Therefore, we are interested in the number of edges in a graph where $\nu(G)$ is bounded by $m$ and $\Delta(G)$ is bounded by $d$.

Triangle-free graphs are defined as graphs with no cycle on three vertices. In other words, the size of the largest complete subgraph in a (non-empty) triangle-free graph is two. Let $\Delta$ denote the class of all triangle-free graphs. Then, the maximum number of edges in a triangle-free graph with maximum degree at most $d$ and maximum matching size at most $m$ is denoted by $f_{\Delta}(d, m)$. 

An $\ell$-star $K_{1,\ell}$ is a bipartite graph with one (central) vertex on one side, which is adjacent to $\ell$ vertices on the other side. If all the vertices of a graph have the same degree $d$, we call it \textit{$d$-regular}. If only one vertex has degree $d-1$ and the rest of the vertices have degree $d$, we call it \textit{almost d-regular}. A graph is called \textit{factor-critical} if removing any vertex from $G$ leaves a graph that admits a perfect matching, that is, a matching saturating all vertices.

The following theorem in \cite{bk09} gives the formula for $f_{GEN}(d, m)$ for all $d$ and $m$ along with a construction for extremal graphs.

\begin{theorem}
\cite{bk09} For general graphs with $\Delta(G) \le d$ and $\nu(G) \le m$, the maximum number of edges in an extremal graph is
$$
f_{GEN}(d, m) = dm + \lfloor{d/2}\rfloor \lfloor{m/\lceil{d/2}\rceil}\rfloor.
$$

Moreover, an extremal graph with $f_{GEN}(d, m)$ edges can be obtained by taking the disjoint union of $r$ copies of $d$-star and $q$ copies of
$$
\begin{cases}
    K_{d+1} & \text{if d+1 is odd} \\
    K^{'}_{d+1} & \text{if d+1 is even},
\end{cases}
$$

\noindent where $q$ is the largest integer such that $m = q \lceil{d/2}\rceil + r$ and $r \ge 0$ and where $K^{'}_{d+1}$
is the graph obtained by removing a perfect matching from the complete graph $K_{d+1}$ on $d + 1$ vertices, adding a new vertex $v$, and making $v$ adjacent to $d$ of the other vertices.
\end{theorem}

The same extremal problem when restricted to triangle-free graphs has been addressed by Ahanjideh, Ekim, and Y\i ld\i z in \cite{aey22}. The authors settle all the cases for $d \ge m$, and for $d < m$ with either $d \le 6$ or $Z(d) \le m < 2d$ where $Z(d)$ is defined as follows.

\begin{definition}\cite{aey22}
For any $d \ge 2$, let $Z(d)$ be the smallest natural number $n$ such that there exists a $d$-regular (if $d$ is even) or almost $d$-regular (if $d$ is odd) triangle-free and
factor-critical graph $G$ with ${\nu(G)}=n$.
\end{definition}

Since $Z(d)$ plays a crucial role in the description of triangle-free extremal graphs, the authors first describe a graph construction that proves the existence of $Z(d)$, then investigate further the value of $Z(d)$.

\begin{lemma}\cite{aey22}
We have $Z(d)=d$ for $d\in \{2,3\}$ and $Z(d)=d+1$ for $d\in\{4,5\}$.
\end{lemma}

\begin{lemma}\cite{aey22}
For $d \ge 2$, if $d$ is even then we have $Z(d) = \lfloor{5d/4}\rfloor$; if $d$ is odd then
we have $\lfloor{5(d-1)/4}\rfloor \le Z(d) \le \lceil{5(d+1)/4}\rceil$.
\label{theorem:z}
\end{lemma}

The main result in \cite{aey22} provides the following formula for $f_{\Delta}(d, m)$ in all solved cases along with a description of extremal graphs, which we omit here for the sake of brevity. 

\begin{theorem}
\cite{aey22} Let $d$ and $m$ be natural numbers with $d \ge 2$, and let $k$ and $r$ be non-negative integers such that $m = kZ(d)+r$ with $0 \le r < Z(d)$. Then, for all cases
with $d \ge m$, and for the cases $d < m$ with either $d \le 6$ or $Z(d) \le m < 2d$, we have

$$
    f_{\Delta}(d, m)= 
\begin{cases}
    dm + k\lfloor{d/2}\rfloor,& \text{if } r < d\\
    dm + k\lfloor{d/2}\rfloor + r - d + 1, & \text{if } r \geq d.
\end{cases}
$$
\label{thm:main}
\end{theorem}

It follows from Theorem \ref{thm:main} that the remaining open cases are $d \ge 7$ with either $m \ge 2d$ or $d < m < Z(d)$. Apart from the formula for $f_{\Delta}(d, m)$, an important contribution in \cite{aey22} is the following result, which shows that there is a triangle-free extremal graph with a special structure. This structural property expressed as a combination of several results (namely Corollary $2.2$ and Lemma $4.4$) in \cite{aey22}, is crucial in building our integer programming formulations. 

\begin{lemma} \cite{aey22}
Let $d$ and $m$ be natural numbers with $d\geq 2$, and let $G$ be an edge-extremal graph with maximum number of connected components isomorphic to a $d$-star. Then, for every connected component $H$ of $G$, one of the following is true:\\
(i) $H$ is a d-star.\\
(ii) $H$ is factor-critical with $|E(H)| = f_\Delta(d, \nu(H))$ and $|V(H)| = 2\nu(H) + 1$ where $d \le \nu(H) \le Z(d)$.
\label{lemma:factor}
\end{lemma}

Now, consider the open case for $m \ge 2d$ and  $d \ge 7$.
Lemma \ref{lemma:factor} states that there is a triangle-free edge-extremal graph whose components which are not $d$-stars are edge-extremal factor-critical triangle-free graphs with matching number between $d$ and $Z(d)$. Let us call the latter {\em extremal components} throughout the paper. This information is used in \cite{aey22} to formulate the open cases with $m \ge 2d$ and  $d \ge 7$ as a knapsack problem where the utility parameters are $f_{\Delta}(d, i)$ for $d\leq i \leq Z(d)$, which are yet to be calculated. Let $x_i$ be the number of extremal components of $G$ with matching number $i$. Then an extremal graph for $d$ and $m$ can be obtained as follows. For $d\leq i \leq Z(d)$, take $x_i$ many extremal components with matching number $i$; their contribution to the total number of edges $f_{\Delta}(d, m)$ is the second summation on the left-hand side of Equation (\ref{eq:sum}), and they contribute $\sum_{i = d}^{Z(d)} ix_i$ to the matching number of the graph. Complete the matching number to $m$ using $d$-stars, each of which adding $d$ to $f_{\Delta}(d,m)$; this is expressed by the first term on the left-hand side of Equation (\ref{eq:sum}). Then we have the following: 
\begin{center}
\begin{equation}
    f_{\Delta}(d, m) = d(m - \sum_{i = d}^{Z(d)} ix_i) + \sum_{i = d}^{Z(d)} f_{\Delta}(d, i) x_i = dm + \sum_{i = d}^{Z(d)} (f_{\Delta}(d, i) - di) x_i. 
\label{eq:sum}
\end{equation}
\end{center}

Based on these observations, Ahanjideh, Ekim, and Y\i ld\i z  propose the following bounded knapsack formulation where the utility of item $i$ is $(f_{\Delta}(d, i) - di)$ and its volume is $i$ \cite{aey22}.

\begin{align}
\label{model:knapsackstart}\text{(Knapsack Formulation)} \quad \max dm + \sum_{i = d}^{Z(d)} &(f_{\Delta}(d, i) - di)  x_i \\
\text{s.t.} \quad \sum_{i = d}^{Z(d)} i x_i &\leq m \\
\label{model:knapsackend} x_{i} \geq 0, \quad x_{i} &\in \mathbb{Z}
\end{align}

 With this formulation, if the edge counts of the extremal components are known for a fixed $d$, then edge extremal graphs for all $m$ can be found. Therefore, in this paper, our effort is focused on finding the extremal components, that is the edge-extremal graphs $H$ for $d < m$ with $d > 6$ and $d < \nu(H) < Z(d)$. Note that the extremal components with matching number $d$ and $Z(d)$ follow from Theorem \ref{thm:main}; we have $f_{\Delta}(d, d)=d^2+1$ and $f_{\Delta}(d, Z(d))=dZ(d)+\lfloor \frac{d}{2} \rfloor$. However, $Z(d)$ is not known when $d$ is odd.

Although the authors in \cite{aey22} leave the development of further methods to solve this knapsack formulation, they still provide some insights about its solution. In particular, they conjecture that its unknown parameters, namely $f_{\Delta}(d, i)$ for $7 \le d<i<Z(d)$ follow the formula given in Theorem \ref{thm:main}, which simplifies as follows in this case:

\begin{conjecture}\cite{aey22}
For $7 \le d<i<Z(d)$, we have $f_{\Delta}(d, i)=di+i-d+1$. 
\label{conj:comp}
\end{conjecture}

The authors also show that if Conjecture \ref{conj:comp} holds then the Knapsack Formulation admits a special optimal solution having as many extremal components as possible with matching number $Z(d)$, and at most one extremal component with smaller matching number.

\begin{proposition}\cite{aey22}
If Conjecture \ref{conj:comp} is true, then for $7 \le d < m < Z(d)$, the Knapsack Formulation admits an optimal solution with $\sum_{i = d}^{Z(d)-1} x_i \le 1$. In other words, $x_{Z(d)}$ is maximized and there is at most one other $x_i$ which is 1 (all the rest being zero). \label{prop:number}
\end{proposition}

Based on the observation that the description of the extremal graphs in Proposition \ref{prop:number} is similar to those provided in Theorem \ref{thm:main} (which we omitted in our paper), the authors also conjecture that the formula in Theorem \ref{thm:main} holds in general (including the open cases).

\begin{conjecture}\cite{aey22}
For all natural numbers $d \ge 2$ and $m$, let $m = kZ(d)+r$. Then we have 

$$
    f_{\Delta}(d, m)= 
\begin{cases}
    dm + k\lfloor{d/2}\rfloor,& \text{if } r < d\\
    dm + k\lfloor{d/2}\rfloor + r - d + 1, & \text{if } r \geq d.
\end{cases}
$$
\label{conj:lemma}
\end{conjecture}

Next, we develop integer programming formulations to construct extremal components. Our method allows us to construct all extremal components for $d=7,8,9,10$ and consequently to solve the Knapsack Formulation for these $d$ values (and any $m$). All our findings support Conjecture \ref{conj:comp}, Proposition \ref{prop:number}, and Conjecture \ref{conj:lemma}; and therefore strengthen them. As a byproduct, we also obtain some new values for $Z(d)$, namely, $Z(7)=9, Z(9)=13, Z(11)=15$, and $Z(13)=17$.

\section{Methodology}

We formulate the construction of an edge extremal triangle-free graph with a matching number at most $m$ and degree at most $d$ as an integer programming problem. Let $V$ be the vertex set and let $x_{ij}$ be a binary variable defined for $i\neq j$ that takes on value 1 if there is an edge between $i, j \in V$ and 0 otherwise. We work with undirected graphs; therefore, we only use one of $x_{ij}$ and $x_{ji}$ variables. For notational simplicity in the models, $ij$ in $x_{ij}$ should be considered an unordered set while the edges are only defined for $i > j$.

\begin{align}
\text{(Basic Formulation)} \quad \max \sum_{i,j \in V} & x_{ij} \label{eq:basic_obj} \\
\text{s.t.} \quad x_{ij} + x_{jk} + x_{ik} &\leq 2  \quad \forall i,j,k \in V \label{eq:triangle} \\
\sum_{j \in V} x_{ij} &\leq d \quad \forall i \in V \label{eq:degree} \\
x_{ij} &\in \{0, 1\} \quad \forall i, j \in V \label{eq:basic_end}
\end{align}

The objective function (\ref{eq:basic_obj}) maximizes the edge count. Constraint (\ref{eq:triangle}) ensures that the resulting graph is triangle-free and Constraint (\ref{eq:degree}) limits the maximum degree. We do not bound the matching number explicitly; instead, we make use of Lemma \ref{lemma:factor}, which guarantees the existence of an edge-extremal triangle-free graph that is factor-critical and fix $|V|=2m+1$ accordingly. Clearly, this bounds the matching number with $m$. Since there is a triangle-free graph with $2m+1$ vertices and degree bounded by $d$ for every $d$ and $m$, an optimal solution of the Basic Formulation describes an extremal triangle-free graph having $f_{\Delta}(d, m)$ edges. The graph resulting from the Basic Formulation might or might not be factor-critical; however, Lemma \ref{lemma:factor} guarantees that it has $f_{\Delta}(d, m)$ edges because there is at least one extremal component which is factor-critical.

The Basic Formulation is highly symmetric since all $x_{ij}$'s are interchangeable. To overcome this problem, we propose two approaches and combine them as a third approach. Our first approach is Orbital Branching, which modifies the branching decisions during the solution of Basic Formulation via the branch-and-cut algorithm. Orbital Branching helps us eliminate symmetric solutions in the branch-and-cut tree. Our second approach is solving the problem with minor modifications iteratively, which we call the Iterative Method. In the first iteration, all the vertex degrees are set to $d$. If there is no feasible solution, all the degrees are set to $d$ but one, which is constrained to be less than $d$, reducing the upper bound of the solution. As long as there is no feasible solution, the degree of one more vertex is constrained to be less than $d$ at every iteration. This approach basically restricts the feasible region of the problem.

\subsection{Orbital Branching}
\label{section_orbital}
We use Orbital Branching to modify the branching decisions while solving the Basic Formulation. The idea is to branch simultaneously on multiple variables on the same symmetry group instead of branching on a single variable at every node of the branch-and-cut tree. 

Before explaining the symmetry groups and equivalency of variables in formal terms, let us first investigate the structure of our Basic Formulation. All $x_{ij}$'s are interchangeable since nothing distinguishes one edge from others in the formulation. Assume that we create a branch for each variable $x_{12} = 1$, $x_{13} = 1$, and so on, as well as a branch where all variables are $0$ to cover the feasible region. This is a legitimate branching strategy since the entire feasible region is covered, even if there are intersections between branches. However, we can observe that each child node with $x_{1j} = 1, j\in V\setminus\{1\}$ yields the same subproblem. Therefore, instead of branching on $|V|$ nodes, we can branch on only two nodes by selecting a representative node, say $x_{12} = 1$, and creating another branch with $x_{12} = x_{13} = ... = 0$. This is the main idea behind Orbital Branching. In the root node, it is easy to see the equivalent variables; but after the first branching, it gets more complicated. Now, $x_{13}$ is not equivalent with $x_{34}$ since they are creating different subproblems due to the already fixed variable $x_{12} = 1$ at that stage. Therefore, we need a formal definition of symmetry groups to identify them. To this end, let us consider the following problem:

\begin{equation}
\max_{x \in \{0,1\}^n} \{c^{T} x | Ax \le b\}. \hspace{3mm}
\end{equation}

A symmetry group of matrix $A$ is defined as the set of permutations which leaves $A$ invariant \cite{m03}. This means that there exists a permutation of variables (columns) and constraints (rows) such that applying them consecutively on $A$ creates the same matrix $A$ \cite{o11}. If the permutation of columns also creates the same $c$ and the permutation of rows creates the same $b$; the new problem is equivalent to the first one. Symmetry groups of an IP can be determined via graph automorphism. The IP is converted into a bipartite graph $B(V, M, E)$, where $V$ represents variables and $M$ represents constraints. There is an edge $(v_i, m_j)$ if only if $A_{ij} \neq 0$, that is, variable $i$ occurs in constraint $j$ with a non-zero coefficient. To account for the differences in $c$, $b$, and $A$ values, a color code is used to show interchangeable vertices in $B$ \cite{pr19}. In general, for every cost coefficient $c_i$ value, a different color is used for vertex $v_i \in V$. Similarly, vertex $m_j \in M$ is colored according to $b_j$ and the inequality type of the constraint. In this setting, changing the order of columns is equivalent to changing the labels of vertices in $V$, and changing the order of rows is equivalent to changing the labels of vertices in $M$. Note that we should prevent interchanging a vertex from $V$ with a vertex from $M$, so they are also assigned different colors. Then, the symmetry groups of the IP and the automorphism groups of the graph $B(V, M, E)$ are the same. If there exists an automorphism assigning $x_i$ to $x_{i'}$, they are deemed equivalent and part of the same orbit. Orbits consist of equivalent variables and they are used for branching decisions in Orbital Branching.

Calculating automorphism groups divides the vertices into orbits, which are set of equivalent vertices. There are two sets of orbits, variable orbits, and constraint orbits. For our Orbital Branching approach, we use the variable orbits and prioritize  the orbit having the maximum number of elements to branch on.

In Figure \ref{fig:branchbound} an illustrative example of Orbital Branching is provided for the model below.

\begin{align}
\max \quad x + y + z + t \label{eq:example_start} \\
\text{s.t.} \quad x + y &\ge 1 \\
x + y &\ge 1 \\
x + y &\ge 1 \\
x, y, z, t &\in \{0,1\} \label{eq:example_end}
\end{align} 

\begin{figure}[h]
    \centering
\begin{tikzpicture}
		\node [style=standard] (0) at (-19, 8) {x};
		\node [style=standard] (1) at (-19, 7) {y};
		\node [style=standard] (2) at (-19, 6) {z};
		\node [style=standard] (3) at (-19, 5) {t};
		\node [style=standard] (4) at (-17, 8) {c1};
		\node [style=standard] (5) at (-17, 7) {c2};
		\node [style=standard] (6) at (-17, 6) {c3};
		\node [style=standard] (7) at (-14, 8) {x};
		\node [style=standard] (8) at (-14, 7) {z};
		\node [style=standard] (9) at (-14, 6) {y};
		\node [style=standard] (10) at (-14, 5) {t};
		\node [style=standard] (11) at (-12, 8) {c2};
		\node [style=standard] (12) at (-12, 7) {c1};
		\node [style=standard] (13) at (-12, 6) {c3};
        \node [style=standard] (17) at (-9, 8) {x};
		\node [style=standard] (18) at (-9, 7) {t};
		\node [style=standard] (19) at (-9, 6) {y};
		\node [style=standard] (20) at (-9, 5) {z};
		\node [style=standard] (21) at (-7, 8) {c3};
		\node [style=standard] (22) at (-7, 7) {c1};
		\node [style=standard] (23) at (-7, 6) {c2};
		\node (14) at (-16.5, 7) {};
		\node (15) at (-14.25, 7) {};
		\node (16) at (-15.5, 7.5) {(y,z)(c1,c2)};
  	\node (24) at (-11.5, 7) {};
		\node (25) at (-9.25, 7) {};
		\node (26) at (-10.5, 7.5) {(z,t)(c2,c3)};
		\draw [style=new edge style 0] (0) to (4);
		\draw [style=new edge style 0] (0) to (5);
		\draw [style=new edge style 0] (0) to (6);
		\draw [style=new edge style 0] (1) to (4);
		\draw [style=new edge style 0] (2) to (5);
		\draw [style=new edge style 0] (3) to (6);
		\draw [style=new edge style 0] (7) to (11);
		\draw [style=new edge style 0] (7) to (12);
		\draw [style=new edge style 0] (7) to (13);
		\draw [style=new edge style 0] (8) to (11);
		\draw [style=new edge style 0] (9) to (12);
		\draw [style=new edge style 0] (10) to (13);
		\draw [style=new edge style 1] (14) to (15);
  	\draw [style=new edge style 0] (17) to (21);
		\draw [style=new edge style 0] (17) to (22);
		\draw [style=new edge style 0] (17) to (23);
		\draw [style=new edge style 0] (18) to (21);
		\draw [style=new edge style 0] (19) to (22);
		\draw [style=new edge style 0] (20) to (23);
		\draw [style=new edge style 1] (24) to (25);
\end{tikzpicture}
\caption{An automorphism of the bipartite graph corresponding to the model (\ref{eq:example_start})-(\ref{eq:example_end})}
\label{fig:automorphism}
\end{figure}
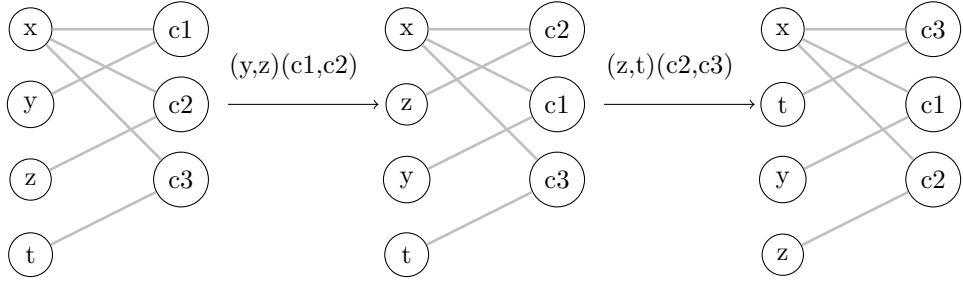

\begin{figure}[h]
    \centering
    \begin{minipage}{0.45\textwidth}
        \centering
        \begin{forest}
          branch and bound,
          where level=1{
            set branch labels={y=z=t=}{}{y=}{},
          }{
            if level=2{
              set branch labels={x_{2,3}=}{}{x_{2,3}=}{},
            }{},
          }
          [1:S:4
            [1:S_1:1:0]
            [2:S_2:4:1]
          ]
        \end{forest}
        \caption{Branch-and-cut tree using Orbital Branching}
        \label{fig:branchbound}
    \end{minipage}\hfill
    \begin{minipage}{0.45\textwidth}
        \centering
        \begin{tikzpicture}
        		\node [style=standard] (0) at (-6, 3) {};
        		\node [style=standard] (1) at (-6, 2) {z};
        		\node [style=standard] (2) at (-6, 1) {t};
        		\node [style=standard] (4) at (-4, 3) {c1};
        		\node [style=standard] (5) at (-4, 2) {c2};
        		\node [style=standard] (6) at (-4, 1) {c3};
        		\node [style=standard] (9) at (-6, 3) {x};
        		\draw [style=new edge style 0] (9) to (4);
        		\draw [style=new edge style 0] (9) to (5);
        		\draw [style=new edge style 0] (9) to (6);
        		\draw [style=new edge style 0] (1) to (5);
        		\draw [style=new edge style 0] (2) to (6);
        \end{tikzpicture}
        \caption{Bipartite graph for $S_2$}
        \label{fig:bipartiteS2}
    \end{minipage}
\end{figure}

In the graph representation of the model (\ref{eq:example_start})-(\ref{eq:example_end}) in Figure \ref{fig:automorphism}, the vertices $x, y, z, t$ represent the variables and vertices $c1, c2, c3$ represent the constraints. Let the graph on the left side of Figure \ref{fig:automorphism} be the original bipartite graph. Then, we can see that interchanging vertex labels $y, z$, and $c1, c2$ create the same incidence matrix. This is equivalent to saying that the graph on the right is obtained by relabeling the vertices of the graph on the left differently, thus the two graphs in Figure \ref{fig:automorphism} are isomorphic. This procedure is represented with $(y,z)(c1,c2)$ and it is an automorphism group of the bipartite graph. Another automorphism group is $(z,t)(c2,c3)$. The application of automorphism groups consecutively yields also an automorphism group. Therefore, we can see that $y, z, t$ are interchangeable and share the same orbit. Consequently, Orbital Branching suggests the branching shown in Figure \ref{fig:branchbound}, where the upper and lower bounds for each node are written over the right and left of each node respectively. In each node of at depth $1$ of the branch-and-cut tree, an optimal solution is found and the orbits are recalculated based on the new graph. The bipartite graph representing the subproblem $S_2$ is given in Figure \ref{fig:bipartiteS2} as an example. Let $F_0^a$ and $F_1^a$ denote the variables set to $0$ and $1$ in node $a$ of a branch-and-cut tree respectively. Then the Orbital Branching algorithm can be summarized as in Algorithm \ref{alg:orbitalDefault}:

\begin{algorithm}
\caption{Orbital Branching}\label{alg:orbitalDefault}
\begin{algorithmic}
\State \textbf{Step 1:} Given a branch-and-cut node $a = (F_0^a, F_1^a)$, calculate the set of orbits $O_a$ at node $a$
\State \textbf{Step 2:} Pick an arbitrary orbit $O \in O_a$
\State \textbf{Step 3:} Pick an arbitrary variable $v_i\in O$, return new nodes $l = (F_0^a \cup \{O\}, F_1^a)$, $r = (F_0^a, F_1^a \cup \{v_i\})$
\end{algorithmic}
\end{algorithm}

We note that all $c$ values are $1$ in our Basic Formulation; therefore we do not need to color variable vertices.

\subsection{Iterative Method}
\label{iterative_section}

Since the maximum degree of a vertex and the number of vertices are limited, the \textit{precomputed upper bound} of the objective function value of the Basic Formulation is $\lfloor{(2m + 1) * d/2)}\rfloor$. This bound can only be achieved if all vertices have the maximum degree $d$. If exactly one of them has degree $d-1$ and all the others have degree $d$, then the upper bound may decrease since the expression value before rounding down decreases by $1/2$. It is possible to use this information to decrease iteratively the upper bound of the Basic  Formulation. For instance, let $d=8$ and $m=9$; thus there are $2*9 + 1 = 19$ vertices. We have $(8*19)/2 = 76$ as an upper bound. We can set all degrees to $8$ and check if there is a feasible solution. If there is a feasible solution, it is an optimal solution for the Basic Formulation. If the problem is infeasible, we know that an upper bound of the problem is $75$ and at least one of the vertices has a degree at most $d-1$. 

Our Iterative Method starts by fixing all the degrees to the degree bound $d$, solving the resulting integer program, and decreasing the upper bound one by one. Each iteration works on a portion of the original feasible region and decreases the upper bound iteratively. The idea is to obtain an optimal solution to the Basic Formulation by solving the formulation in a restricted feasible region. We can consider different formulations of this approach, enumerating all possible degree combinations for each objective value or simply decreasing the upper bound by setting the degrees of some vertices to be less than $d$ while the rest are equal to $d$. We adopt the latter formulation since enumerating all possible degrees would imply an exponential number of iterations.

We distinguish two sets of vertices; let $V_d$ be the set of vertices of degree equal to $d$ and $V_{d-1}$ be the set of vertices of degree at most $d-1$. The formulation can be given as:
\begin{align}
\text{(Iterative Formulation)} \quad \max \sum_{i,j \in V} & x_{ij} \label{it_start} \\
\text{s.t.} \quad x_{ij} + x_{jk} + x_{ik} &\leq 2  \quad \forall i,j,k \in V \\
\sum_{j \in V} x_{ij} &= d \quad \forall i \in V_d, \label{constVm} \\
\sum_{j \in V} x_{ij} &\leq d-1 \quad \forall i \in V_{d-1}, \label{constVn} \\
x_{ij} &\in \{0, 1\} \quad \forall i, j \in V \label{it_end}
\end{align} 

It should be noted that an optimal solution to the Iterative Formulation with $|V_{d-1}|=1$ might not be an optimal solution for the Basic Formulation. Indeed, a graph where all vertices have degree $d$ but one which has degree $d-4$ has fewer edges than a graph with $|V(G)|-2$ vertices of degree $d$, and two vertices of degree $d-1$. Consequently, the optimal solution of the Iterative Formulation with $|V_{d-1}|=1$ gives a lower bound, which can possibly be improved later by adding more vertices in $V_{d-1}$. In each iteration, sets $V_d$ and $V_{d-1}$ are updated, and the Iterative Formulation is solved. The Iterative Method summarized in Algorithm \ref{alg:rapidIteration} solves the Iterative Formulation successively by updating the upper bounds until the upper bound is equal to the lower bound obtained by the best feasible solution. Let $N=2m+1$ be the number of vertices in what follows. 

\begin{algorithm}
\caption{Iterative Method}\label{alg:rapidIteration}
\begin{algorithmic}
\State \textbf{Input}:  $V_{d-1} = \emptyset$, $V_d = V$
\State $LB = 0$, $UB =  \lfloor{((Nd) / 2} \rfloor$
\While{$UB > LB$}
\State Solve Iterative Formulation $D$ with $V_d$ and $V_{d-1}$. 
\If{$D$ is feasible}
\State $LB = \max(z^{*}, LB)$ where $z^{*}$ is the optimal objective function value of $D$
\EndIf
\State Pick $u \in V_d$
\State $V_d = V_d \setminus \{u\}$
\State $V_{d-1} = V_{d-1} \cup \{u\}$
\State $UB = \max (\lfloor{((Nd) / 2 - 0.5|V_{d-1}|}\rfloor$, $LB$)
\EndWhile
\State \Return $UB$ and the optimal solution of the Iterative Formulation giving the best $LB$.
\end{algorithmic}
\end{algorithm}

Our Iterative Method does not primarily focus on reducing symmetry, instead, it focuses on reducing the size of the feasible region. In the first iteration, where $V_{d-1} = \emptyset$, all variables are equivalent. Still, in each iteration, one of the degree constraints for a vertex changes, reducing the symmetry. To reduce the symmetry further, we can order the vertices of $V_{d-1}$ lexicographically so that their degrees are non-increasing.

\subsection{Iterative Method with Orbital Branching}
\label{combined_section}
Orbital Branching exploits the symmetry for branching decisions, and the Iterative Method reduces the size of the feasible region while keeping most of the symmetry intact. Besides, the Orbital Branching interferes with the default branch-and-cut decisions of the solver while the Iterative Method modifies the problem formulation. Therefore it is possible to use these methods jointly so that the same iterations as in the Iterative Method are done, but branching decisions are made with Orbital Branching. The Iterative Formulation is less symmetric than the Basic Formulation since Equations (\ref{constVm}) are equalities. 

\section{Computational Experiments}

Computations have been conducted on Intel(R) Core(TM) i7-9750H CPU @ 2.60GHz with 32 GB RAM using CPLEX 20.1.0 on $10$ threads. We used CPLEX 20.1.0 with C++ to implement the methods. We enforce a limit of $1800$ seconds for each run. In the following tables that show our computational results, ``PreUB" stands for the precomputed upper bound, ``LB" stands for the lower bound found by CPLEX, ``UB" stands for the upper bound found by CPLEX, ``Gap" represents the relative optimality gap, ``Time" is the execution time, and ``Node" is the number of nodes in the branch-and-cut tree. In Tables \ref{table:basiccp},\ref{table:orbital}, and \ref{table:iterperf}, lines corresponding to optimality are marked with a dark background. We tested instances with $7 \le d \le 13$ and $d < m \le Z(d)$. When $d$ is odd, $Z(d)$ is not known beforehand. However, by definition of $Z(d)$, a graph with matching number $m=Z(d)$ is almost $d$-regular. Accordingly, we find $Z(d)$ based on the degrees in the resulting edge-extremal graph. More precisely, $Z(d)$ is the value of the matching number $m$ for which the LB is equal to the precomputed UB and the resulting extremal graph is factor-critical, which we check by computation. When $d$ is even, we have $Z(d)=\lfloor 5d/4\rfloor$ by Lemma \ref{theorem:z}. This implies $Z(8)=10, Z(10)=12$ and $Z(12)=15$. Then, it follows from Theorem \ref{thm:main} that $f_{\Delta}(8, 10)=84$, $f_{\Delta}(10, 12)=125$ and $f_{\Delta}(12, 15)=186$. From the computations, we learned the $Z(d)$ values for the new $d$ values 7, 9, 11, and 13; namely $Z(7) = 9$, $Z(9) = 12$, $Z(11) = 15$, $Z(13) = 17$. Overall, if we check the solved instances, the edge counts align with the Conjecture \ref{conj:lemma} given in Ahanjideh, Ekim, and Yıldız \cite{aey22}.

For solving the Basic Formulation (\ref{eq:basic_obj})-(\ref{eq:basic_end}), we started with CPLEX with default parameters and CPLEX with symmetry breaking parameter set to its most aggressive setting (\textit{IloCplex::Param::Preprocessing::Symmetry} $=5$). The results are shared in Table \ref{table:basiccp}. Both of the approaches solved only 7 of the instances to optimality. Overall, CPLEX with symmetry breaking is slightly better due to a lower average time. Gaps are equivalent for both of the approaches. Note that the upper bound levels are equivalent to precomputed upper bound values that can be achieved by the maximum degree constraint. For example, when $d$ is $7$ and $m$ is $8$, there are $17$ vertices since we can assume that an extremal graph is factor-critical. Since $d$ is odd, the maximum edge count is achieved via an almost $d$-regular graph (the sum of all the degrees should be even). $16$ vertices with degree $7$ and one vertex with degree $6$, yielding $59$ edges in total. This explains why the only solved instances are for $m = Z(d)$ values. When $m = Z(d)$, the resulting graph is $d$-regular or almost $d$-regular, which has an equivalent edge count to the bound we discussed. From an inspection of the results, we can conclude that our initial approaches are not good enough for decreasing the upper bound found by CPLEX, but they are promising for finding feasible solutions. This might be due to symmetry. In our next experiment, we test the efficacy of using Orbital Branching with Basic Formulation (Subsection \ref{section_orbital}).

\begin{table}[h]
\centering
\resizebox{\textwidth}{!}{
\begin{tabular}{|ccc|ccccc|ccccc|}
\hline
\multicolumn{3}{|c|}{\textbf{Parameters}}                                                                           & \multicolumn{5}{c|}{\textbf{Basic Formulation}}                                                                                                                                                                               & \multicolumn{5}{c|}{\textbf{Basic Formulation + Symmetry Breaking}}                                                                                                                                                          \\ \hline
\multicolumn{1}{|c|}{\textit{\textbf{d}}}        & \multicolumn{1}{c|}{\textit{\textbf{m}}}        & \textbf{PreUB} & \multicolumn{1}{c|}{\textbf{LB}}                 & \multicolumn{1}{c|}{\textbf{UB}}                 & \multicolumn{1}{c|}{\textbf{Gap}}                   & \multicolumn{1}{c|}{\textbf{Time}}                & \textbf{Node} & \multicolumn{1}{c|}{\textbf{LB}}                 & \multicolumn{1}{c|}{\textbf{UB}}                 & \multicolumn{1}{c|}{\textbf{Gap}}                   & \multicolumn{1}{c|}{\textbf{Time}}               & \textbf{Node} \\ \hline
\multicolumn{1}{|c|}{7}                          & \multicolumn{1}{c|}{8}                          & 59             & \multicolumn{1}{c|}{58}                          & \multicolumn{1}{c|}{59}                          & \multicolumn{1}{c|}{1.72\%}                         & \multicolumn{1}{c|}{1800}                         & 8412700       & \multicolumn{1}{c|}{58}                          & \multicolumn{1}{c|}{59}                          & \multicolumn{1}{c|}{1.72\%}                         & \multicolumn{1}{c|}{1800}                        & 4657959       \\ \hline
\rowcolor[HTML]{C0C0C0} 
\multicolumn{1}{|c|}{\cellcolor[HTML]{C0C0C0}7}  & \multicolumn{1}{c|}{\cellcolor[HTML]{C0C0C0}9}  & 66             & \multicolumn{1}{c|}{\cellcolor[HTML]{C0C0C0}66}  & \multicolumn{1}{c|}{\cellcolor[HTML]{C0C0C0}66}  & \multicolumn{1}{c|}{\cellcolor[HTML]{C0C0C0}0.00\%} & \multicolumn{1}{c|}{\cellcolor[HTML]{C0C0C0}0}    & 0             & \multicolumn{1}{c|}{\cellcolor[HTML]{C0C0C0}66}  & \multicolumn{1}{c|}{\cellcolor[HTML]{C0C0C0}66}  & \multicolumn{1}{c|}{\cellcolor[HTML]{C0C0C0}0.00\%} & \multicolumn{1}{c|}{\cellcolor[HTML]{C0C0C0}0}   & 0             \\ \hline
\multicolumn{1}{|c|}{8}                          & \multicolumn{1}{c|}{9}                          & 76             & \multicolumn{1}{c|}{74}                          & \multicolumn{1}{c|}{76}                          & \multicolumn{1}{c|}{2.70\%}                         & \multicolumn{1}{c|}{1800}                         & 4605513       & \multicolumn{1}{c|}{74}                          & \multicolumn{1}{c|}{76}                          & \multicolumn{1}{c|}{2.70\%}                         & \multicolumn{1}{c|}{1800}                        & 3012473       \\ \hline
\rowcolor[HTML]{C0C0C0} 
\multicolumn{1}{|c|}{\cellcolor[HTML]{C0C0C0}8}  & \multicolumn{1}{c|}{\cellcolor[HTML]{C0C0C0}10} & 84             & \multicolumn{1}{c|}{\cellcolor[HTML]{C0C0C0}84}  & \multicolumn{1}{c|}{\cellcolor[HTML]{C0C0C0}84}  & \multicolumn{1}{c|}{\cellcolor[HTML]{C0C0C0}0.00\%} & \multicolumn{1}{c|}{\cellcolor[HTML]{C0C0C0}0}    & 0             & \multicolumn{1}{c|}{\cellcolor[HTML]{C0C0C0}84}  & \multicolumn{1}{c|}{\cellcolor[HTML]{C0C0C0}84}  & \multicolumn{1}{c|}{\cellcolor[HTML]{C0C0C0}0.00\%} & \multicolumn{1}{c|}{\cellcolor[HTML]{C0C0C0}0}   & 0             \\ \hline
\multicolumn{1}{|c|}{9}                          & \multicolumn{1}{c|}{10}                         & 94             & \multicolumn{1}{c|}{92}                          & \multicolumn{1}{c|}{94}                          & \multicolumn{1}{c|}{2.17\%}                         & \multicolumn{1}{c|}{1800}                         & 2952258       & \multicolumn{1}{c|}{92}                          & \multicolumn{1}{c|}{94}                          & \multicolumn{1}{c|}{2.17\%}                         & \multicolumn{1}{c|}{1800}                        & 1836200       \\ \hline
\multicolumn{1}{|c|}{9}                          & \multicolumn{1}{c|}{11}                         & 103            & \multicolumn{1}{c|}{102}                         & \multicolumn{1}{c|}{103}                         & \multicolumn{1}{c|}{0.98\%}                         & \multicolumn{1}{c|}{1800}                         & 1356502       & \multicolumn{1}{c|}{102}                         & \multicolumn{1}{c|}{103}                         & \multicolumn{1}{c|}{0.98\%}                         & \multicolumn{1}{c|}{1800}                        & 1122091       \\ \hline
\rowcolor[HTML]{C0C0C0} 
\multicolumn{1}{|c|}{\cellcolor[HTML]{C0C0C0}9}  & \multicolumn{1}{c|}{\cellcolor[HTML]{C0C0C0}12} & 112            & \multicolumn{1}{c|}{\cellcolor[HTML]{C0C0C0}112} & \multicolumn{1}{c|}{\cellcolor[HTML]{C0C0C0}112} & \multicolumn{1}{c|}{\cellcolor[HTML]{C0C0C0}0.00\%} & \multicolumn{1}{c|}{\cellcolor[HTML]{C0C0C0}2}    & 428           & \multicolumn{1}{c|}{\cellcolor[HTML]{C0C0C0}112} & \multicolumn{1}{c|}{\cellcolor[HTML]{C0C0C0}112} & \multicolumn{1}{c|}{\cellcolor[HTML]{C0C0C0}0.00\%} & \multicolumn{1}{c|}{\cellcolor[HTML]{C0C0C0}3}   & 465           \\ \hline
\multicolumn{1}{|c|}{10}                         & \multicolumn{1}{c|}{11}                         & 115            & \multicolumn{1}{c|}{112}                         & \multicolumn{1}{c|}{115}                         & \multicolumn{1}{c|}{2.68\%}                         & \multicolumn{1}{c|}{1800}                         & 2049077       & \multicolumn{1}{c|}{112}                         & \multicolumn{1}{c|}{115}                         & \multicolumn{1}{c|}{2.68\%}                         & \multicolumn{1}{c|}{1800}                        & 1274473       \\ \hline
\rowcolor[HTML]{C0C0C0} 
\multicolumn{1}{|c|}{\cellcolor[HTML]{C0C0C0}10} & \multicolumn{1}{c|}{\cellcolor[HTML]{C0C0C0}12} & 125            & \multicolumn{1}{c|}{\cellcolor[HTML]{C0C0C0}125} & \multicolumn{1}{c|}{\cellcolor[HTML]{C0C0C0}125} & \multicolumn{1}{c|}{\cellcolor[HTML]{C0C0C0}0.00\%} & \multicolumn{1}{c|}{\cellcolor[HTML]{C0C0C0}0}    & 0             & \multicolumn{1}{c|}{\cellcolor[HTML]{C0C0C0}125} & \multicolumn{1}{c|}{\cellcolor[HTML]{C0C0C0}125} & \multicolumn{1}{c|}{\cellcolor[HTML]{C0C0C0}0.00\%} & \multicolumn{1}{c|}{\cellcolor[HTML]{C0C0C0}0}   & 0             \\ \hline
\multicolumn{1}{|c|}{11}                         & \multicolumn{1}{c|}{12}                         & 137            & \multicolumn{1}{c|}{134}                         & \multicolumn{1}{c|}{137}                         & \multicolumn{1}{c|}{2.24\%}                         & \multicolumn{1}{c|}{1800}                         & 530890        & \multicolumn{1}{c|}{134}                         & \multicolumn{1}{c|}{137}                         & \multicolumn{1}{c|}{2.24\%}                         & \multicolumn{1}{c|}{1800}                        & 809227        \\ \hline
\multicolumn{1}{|c|}{11}                         & \multicolumn{1}{c|}{13}                         & 148            & \multicolumn{1}{c|}{146}                         & \multicolumn{1}{c|}{148}                         & \multicolumn{1}{c|}{1.37\%}                         & \multicolumn{1}{c|}{1800}                         & 293767        & \multicolumn{1}{c|}{146}                         & \multicolumn{1}{c|}{148}                         & \multicolumn{1}{c|}{1.37\%}                         & \multicolumn{1}{c|}{1800}                        & 311019        \\ \hline
\multicolumn{1}{|c|}{11}                         & \multicolumn{1}{c|}{14}                         & 159            & \multicolumn{1}{c|}{158}                         & \multicolumn{1}{c|}{159}                         & \multicolumn{1}{c|}{0.63\%}                         & \multicolumn{1}{c|}{1800}                         & 257540        & \multicolumn{1}{c|}{158}                         & \multicolumn{1}{c|}{159}                         & \multicolumn{1}{c|}{0.63\%}                         & \multicolumn{1}{c|}{1800}                        & 257540        \\ \hline
\rowcolor[HTML]{C0C0C0} 
\multicolumn{1}{|c|}{\cellcolor[HTML]{C0C0C0}11} & \multicolumn{1}{c|}{\cellcolor[HTML]{C0C0C0}15} & 170            & \multicolumn{1}{c|}{\cellcolor[HTML]{C0C0C0}170} & \multicolumn{1}{c|}{\cellcolor[HTML]{C0C0C0}170} & \multicolumn{1}{c|}{\cellcolor[HTML]{C0C0C0}0.00\%} & \multicolumn{1}{c|}{\cellcolor[HTML]{C0C0C0}0}    & 0             & \multicolumn{1}{c|}{\cellcolor[HTML]{C0C0C0}170} & \multicolumn{1}{c|}{\cellcolor[HTML]{C0C0C0}170} & \multicolumn{1}{c|}{\cellcolor[HTML]{C0C0C0}0.00\%} & \multicolumn{1}{c|}{\cellcolor[HTML]{C0C0C0}0}   & 0             \\ \hline
\multicolumn{1}{|c|}{12}                         & \multicolumn{1}{c|}{13}                         & 162            & \multicolumn{1}{c|}{158}                         & \multicolumn{1}{c|}{162}                         & \multicolumn{1}{c|}{2.53\%}                         & \multicolumn{1}{c|}{1800}                         & 471193        & \multicolumn{1}{c|}{158}                         & \multicolumn{1}{c|}{162}                         & \multicolumn{1}{c|}{2.53\%}                         & \multicolumn{1}{c|}{1800}                        & 392793        \\ \hline
\multicolumn{1}{|c|}{12}                         & \multicolumn{1}{c|}{14}                         & 174            & \multicolumn{1}{c|}{171}                         & \multicolumn{1}{c|}{174}                         & \multicolumn{1}{c|}{1.75\%}                         & \multicolumn{1}{c|}{1800}                         & 271165        & \multicolumn{1}{c|}{171}                         & \multicolumn{1}{c|}{174}                         & \multicolumn{1}{c|}{1.75\%}                         & \multicolumn{1}{c|}{1800}                        & 218547        \\ \hline
\rowcolor[HTML]{C0C0C0} 
\multicolumn{1}{|c|}{\cellcolor[HTML]{C0C0C0}12} & \multicolumn{1}{c|}{\cellcolor[HTML]{C0C0C0}15} & 186            & \multicolumn{1}{c|}{\cellcolor[HTML]{C0C0C0}186} & \multicolumn{1}{c|}{\cellcolor[HTML]{C0C0C0}186} & \multicolumn{1}{c|}{\cellcolor[HTML]{C0C0C0}0.00\%} & \multicolumn{1}{c|}{\cellcolor[HTML]{C0C0C0}53}   & 6077          & \multicolumn{1}{c|}{\cellcolor[HTML]{C0C0C0}186} & \multicolumn{1}{c|}{\cellcolor[HTML]{C0C0C0}186} & \multicolumn{1}{c|}{\cellcolor[HTML]{C0C0C0}0.00\%} & \multicolumn{1}{c|}{\cellcolor[HTML]{C0C0C0}707} & 82348         \\ \hline
\multicolumn{1}{|c|}{13}                         & \multicolumn{1}{c|}{14}                         & 188            & \multicolumn{1}{c|}{184}                         & \multicolumn{1}{c|}{188}                         & \multicolumn{1}{c|}{2.17\%}                         & \multicolumn{1}{c|}{1800}                         & 326937        & \multicolumn{1}{c|}{184}                         & \multicolumn{1}{c|}{188}                         & \multicolumn{1}{c|}{2.17\%}                         & \multicolumn{1}{c|}{1800}                        & 292989        \\ \hline
\multicolumn{1}{|c|}{13}                         & \multicolumn{1}{c|}{15}                         & 201            & \multicolumn{1}{c|}{198}                         & \multicolumn{1}{c|}{201}                         & \multicolumn{1}{c|}{1.52\%}                         & \multicolumn{1}{c|}{1800}                         & 265385        & \multicolumn{1}{c|}{198}                         & \multicolumn{1}{c|}{201}                         & \multicolumn{1}{c|}{1.52\%}                         & \multicolumn{1}{c|}{1800}                        & 173732        \\ \hline
\multicolumn{1}{|c|}{13}                         & \multicolumn{1}{c|}{16}                         & 214            & \multicolumn{1}{c|}{212}                         & \multicolumn{1}{c|}{214}                         & \multicolumn{1}{c|}{0.94\%}                         & \multicolumn{1}{c|}{1800}                         & 148400        & \multicolumn{1}{c|}{212}                         & \multicolumn{1}{c|}{214}                         & \multicolumn{1}{c|}{0.94\%}                         & \multicolumn{1}{c|}{1800}                        & 144746        \\ \hline
\rowcolor[HTML]{C0C0C0} 
\multicolumn{1}{|c|}{\cellcolor[HTML]{C0C0C0}13} & \multicolumn{1}{c|}{\cellcolor[HTML]{C0C0C0}17} & 227            & \multicolumn{1}{c|}{\cellcolor[HTML]{C0C0C0}227} & \multicolumn{1}{c|}{\cellcolor[HTML]{C0C0C0}227} & \multicolumn{1}{c|}{\cellcolor[HTML]{C0C0C0}0.00\%} & \multicolumn{1}{c|}{\cellcolor[HTML]{C0C0C0}1200} & 45158         & \multicolumn{1}{c|}{\cellcolor[HTML]{C0C0C0}227} & \multicolumn{1}{c|}{\cellcolor[HTML]{C0C0C0}227} & \multicolumn{1}{c|}{\cellcolor[HTML]{C0C0C0}0.00\%} & \multicolumn{1}{c|}{\cellcolor[HTML]{C0C0C0}21}  & 311           \\ \hline
\multicolumn{3}{|c|}{\textbf{Avg}}                                                                                  & \multicolumn{1}{c|}{138.45}                      & \multicolumn{1}{c|}{140}                         & \multicolumn{1}{c|}{1.12\%}                         & \multicolumn{1}{c|}{1233}                         & 1099650       & \multicolumn{1}{c|}{138.45}                      & \multicolumn{1}{c|}{140}                         & \multicolumn{1}{c|}{1.12\%}                         & \multicolumn{1}{c|}{1207}                        & 729346        \\ \hline
\end{tabular}}
\caption{Basic formulation performance summary.}
\label{table:basiccp}
\end{table}

\begin{figure}[H]
\centering
    \begin{minipage}{0.45\textwidth}
    \includegraphics[scale=0.7]{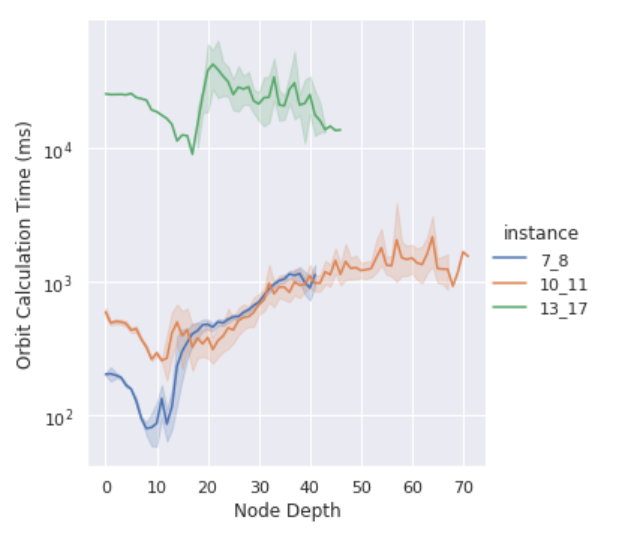}
    \caption{nauty runtime with depth.}
    \label{fig:orbitDuration}
    \end{minipage}\hfill
    \begin{minipage}{0.45\textwidth}
    \includegraphics[scale=0.7]{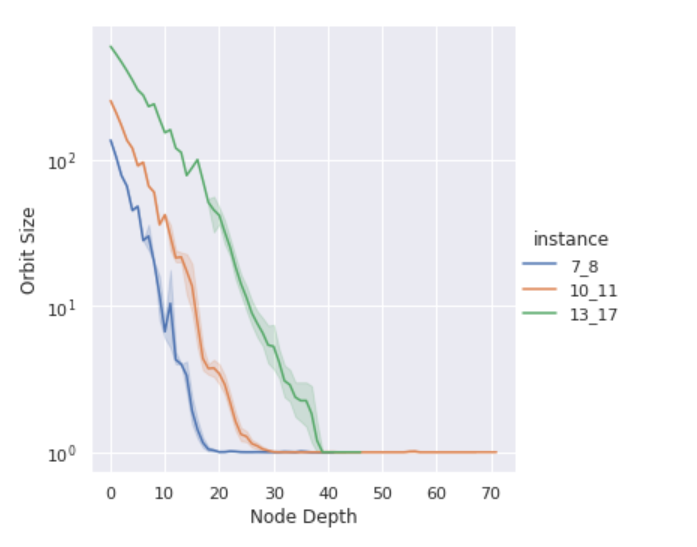}
    \caption{Orbit sizes with depth.}
    \label{fig:orbitSize}
    \end{minipage}\hfill
\end{figure}

We implemented Orbital Branching using the callback mechanism of CPLEX. For finding orbits, we used nauty 2.7r3, a software library for computing automorphism groups of graphs by McKay and Piperno \cite{mp14}. At each branch-and-cut node, we find the variables set to $0$ and $1$ from CPLEX and update the constraints accordingly. Then, we construct the updated bipartite graph using binding constraints. The resulting graph (and coloring) is sent to nauty. For calculating orbits, nauty calculates the automorphism groups utilizing graph isomorphism problem, which is a $GI-complete$ problem \cite{j05}. This implies that the problem of calculating orbits is neither known to be NP-complete nor known to be polynomial-time solvable. For large instances, run time can be as high as $150$ seconds, while it is around $0.1$ seconds for small instances. This can be seen in Figure \ref{fig:orbitDuration}, where instance $7\_8$ means instance with $d = 7$ and $m = 8$. Since there can be tens of thousands of nodes in a branch-and-cut tree, even a runtime of one second for orbit calculation results in an excessive CPU time. Therefore, we analyze the average orbit size for each depth in a branch-and-cut tree. Orbits tend to get smaller as we go deeper in the branch-and-cut tree, which can be seen from Figure \ref{fig:orbitSize}. For example, using Orbital Branching after depth $20$ only adds overhead to the instance with $d = 7$ and $m=8$. Due to this, we run Orbital Branching until a certain depth. This depth is initialized as a large number in the beginning. When we find a node with orbit size $1$, the depth of that node becomes the limit; for deeper nodes, orbits are not calculated. After that depth, the default CPLEX branching strategy is used. Branching is done on the orbit with the maximum number of vertices since it sets the values for more variables in the child nodes. If the largest orbit only has a single variable, we use the default CPLEX branching strategy.

\begin{table}[htb]
\centering
\begin{tabular}{|ccc|ccccc|}
\hline
\multicolumn{3}{|c|}{\textbf{Parameters}}                                                                                    & \multicolumn{5}{c|}{\textbf{Orbital Branching}}                                                                                                                                                                              \\ \hline
\multicolumn{1}{|c|}{\textit{\textbf{d}}}        & \multicolumn{1}{c|}{\textit{\textbf{m}}}        & \textbf{PreUB} & \multicolumn{1}{c|}{\textbf{LB}}                 & \multicolumn{1}{c|}{\textbf{UB}}                 & \multicolumn{1}{c|}{\textbf{Gap}}                   & \multicolumn{1}{c|}{\textbf{Time}}               & \textbf{Node} \\ \hline
\rowcolor[HTML]{C0C0C0} 
\multicolumn{1}{|c|}{\cellcolor[HTML]{C0C0C0}7}  & \multicolumn{1}{c|}{\cellcolor[HTML]{C0C0C0}8}  & 59                      & \multicolumn{1}{c|}{\cellcolor[HTML]{C0C0C0}58}  & \multicolumn{1}{c|}{\cellcolor[HTML]{C0C0C0}58}  & \multicolumn{1}{c|}{\cellcolor[HTML]{C0C0C0}0.00\%} & \multicolumn{1}{c|}{\cellcolor[HTML]{C0C0C0}14}  & 57945         \\ \hline
\rowcolor[HTML]{C0C0C0} 
\multicolumn{1}{|c|}{\cellcolor[HTML]{C0C0C0}7}  & \multicolumn{1}{c|}{\cellcolor[HTML]{C0C0C0}9}  & 66                      & \multicolumn{1}{c|}{\cellcolor[HTML]{C0C0C0}66}  & \multicolumn{1}{c|}{\cellcolor[HTML]{C0C0C0}66}  & \multicolumn{1}{c|}{\cellcolor[HTML]{C0C0C0}0.00\%} & \multicolumn{1}{c|}{\cellcolor[HTML]{C0C0C0}0}   & 0             \\ \hline
\rowcolor[HTML]{C0C0C0} 
\multicolumn{1}{|c|}{\cellcolor[HTML]{C0C0C0}8}  & \multicolumn{1}{c|}{\cellcolor[HTML]{C0C0C0}9}  & 76                      & \multicolumn{1}{c|}{\cellcolor[HTML]{C0C0C0}74}  & \multicolumn{1}{c|}{\cellcolor[HTML]{C0C0C0}74}  & \multicolumn{1}{c|}{\cellcolor[HTML]{C0C0C0}0.00\%} & \multicolumn{1}{c|}{\cellcolor[HTML]{C0C0C0}61}  & 236040        \\ \hline
\rowcolor[HTML]{C0C0C0} 
\multicolumn{1}{|c|}{\cellcolor[HTML]{C0C0C0}8}  & \multicolumn{1}{c|}{\cellcolor[HTML]{C0C0C0}10} & 84                      & \multicolumn{1}{c|}{\cellcolor[HTML]{C0C0C0}84}  & \multicolumn{1}{c|}{\cellcolor[HTML]{C0C0C0}84}  & \multicolumn{1}{c|}{\cellcolor[HTML]{C0C0C0}0.00\%} & \multicolumn{1}{c|}{\cellcolor[HTML]{C0C0C0}0}   & 0             \\ \hline
\rowcolor[HTML]{C0C0C0} 
\multicolumn{1}{|c|}{\cellcolor[HTML]{C0C0C0}9}  & \multicolumn{1}{c|}{\cellcolor[HTML]{C0C0C0}10} & 94                      & \multicolumn{1}{c|}{\cellcolor[HTML]{C0C0C0}92}  & \multicolumn{1}{c|}{\cellcolor[HTML]{C0C0C0}92}  & \multicolumn{1}{c|}{\cellcolor[HTML]{C0C0C0}0.00\%} & \multicolumn{1}{c|}{\cellcolor[HTML]{C0C0C0}558} & 1781858       \\ \hline
\multicolumn{1}{|c|}{9}                          & \multicolumn{1}{c|}{11}                         & 103                     & \multicolumn{1}{c|}{102}                         & \multicolumn{1}{c|}{103}                         & \multicolumn{1}{c|}{0.98\%}                         & \multicolumn{1}{c|}{1800}                        & 2651099       \\ \hline
\rowcolor[HTML]{C0C0C0} 
\multicolumn{1}{|c|}{\cellcolor[HTML]{C0C0C0}9}  & \multicolumn{1}{c|}{\cellcolor[HTML]{C0C0C0}12} & 112                     & \multicolumn{1}{c|}{\cellcolor[HTML]{C0C0C0}112} & \multicolumn{1}{c|}{\cellcolor[HTML]{C0C0C0}112} & \multicolumn{1}{c|}{\cellcolor[HTML]{C0C0C0}0.00\%} & \multicolumn{1}{c|}{\cellcolor[HTML]{C0C0C0}21}  & 193           \\ \hline
\multicolumn{1}{|c|}{10}                         & \multicolumn{1}{c|}{11}                         & 115                     & \multicolumn{1}{c|}{112}                         & \multicolumn{1}{c|}{114}                         & \multicolumn{1}{c|}{1.79\%}                         & \multicolumn{1}{c|}{1800}                        & 1838556       \\ \hline
\rowcolor[HTML]{C0C0C0} 
\multicolumn{1}{|c|}{\cellcolor[HTML]{C0C0C0}10} & \multicolumn{1}{c|}{\cellcolor[HTML]{C0C0C0}12} & 125                     & \multicolumn{1}{c|}{\cellcolor[HTML]{C0C0C0}125} & \multicolumn{1}{c|}{\cellcolor[HTML]{C0C0C0}125} & \multicolumn{1}{c|}{\cellcolor[HTML]{C0C0C0}0.00\%} & \multicolumn{1}{c|}{\cellcolor[HTML]{C0C0C0}0}   & 0             \\ \hline
\multicolumn{1}{|c|}{11}                         & \multicolumn{1}{c|}{12}                         & 137                     & \multicolumn{1}{c|}{134}                         & \multicolumn{1}{c|}{137}                         & \multicolumn{1}{c|}{2.24\%}                         & \multicolumn{1}{c|}{1800}                        & 1223806       \\ \hline
\multicolumn{1}{|c|}{11}                         & \multicolumn{1}{c|}{13}                         & 148                     & \multicolumn{1}{c|}{146}                         & \multicolumn{1}{c|}{148}                         & \multicolumn{1}{c|}{1.37\%}                         & \multicolumn{1}{c|}{1800}                        & 1142867       \\ \hline
\multicolumn{1}{|c|}{11}                         & \multicolumn{1}{c|}{14}                         & 159                     & \multicolumn{1}{c|}{158}                         & \multicolumn{1}{c|}{159}                         & \multicolumn{1}{c|}{0.63\%}                         & \multicolumn{1}{c|}{1800}                        & 428121        \\ \hline
\rowcolor[HTML]{C0C0C0} 
\multicolumn{1}{|c|}{\cellcolor[HTML]{C0C0C0}11} & \multicolumn{1}{c|}{\cellcolor[HTML]{C0C0C0}15} & 170                     & \multicolumn{1}{c|}{\cellcolor[HTML]{C0C0C0}170} & \multicolumn{1}{c|}{\cellcolor[HTML]{C0C0C0}170} & \multicolumn{1}{c|}{\cellcolor[HTML]{C0C0C0}0.00\%} & \multicolumn{1}{c|}{\cellcolor[HTML]{C0C0C0}0}   & 0             \\ \hline
\multicolumn{1}{|c|}{12}                         & \multicolumn{1}{c|}{13}                         & 162                     & \multicolumn{1}{c|}{158}                         & \multicolumn{1}{c|}{162}                         & \multicolumn{1}{c|}{2.53\%}                         & \multicolumn{1}{c|}{1800}                        & 567720        \\ \hline
\multicolumn{1}{|c|}{12}                         & \multicolumn{1}{c|}{14}                         & 174                     & \multicolumn{1}{c|}{171}                         & \multicolumn{1}{c|}{174}                         & \multicolumn{1}{c|}{1.75\%}                         & \multicolumn{1}{c|}{1800}                        & 444420        \\ \hline
\rowcolor[HTML]{C0C0C0} 
\multicolumn{1}{|c|}{\cellcolor[HTML]{C0C0C0}12} & \multicolumn{1}{c|}{\cellcolor[HTML]{C0C0C0}15} & 186                     & \multicolumn{1}{c|}{\cellcolor[HTML]{C0C0C0}186} & \multicolumn{1}{c|}{\cellcolor[HTML]{C0C0C0}186} & \multicolumn{1}{c|}{\cellcolor[HTML]{C0C0C0}0.00\%} & \multicolumn{1}{c|}{\cellcolor[HTML]{C0C0C0}361} & 4769          \\ \hline
\multicolumn{1}{|c|}{13}                         & \multicolumn{1}{c|}{14}                         & 188                     & \multicolumn{1}{c|}{184}                         & \multicolumn{1}{c|}{188}                         & \multicolumn{1}{c|}{2.17\%}                         & \multicolumn{1}{c|}{1800}                        & 278258        \\ \hline
\multicolumn{1}{|c|}{13}                         & \multicolumn{1}{c|}{15}                         & 201                     & \multicolumn{1}{c|}{198}                         & \multicolumn{1}{c|}{201}                         & \multicolumn{1}{c|}{1.52\%}                         & \multicolumn{1}{c|}{1800}                        & 85467         \\ \hline
\multicolumn{1}{|c|}{13}                         & \multicolumn{1}{c|}{16}                         & 214                     & \multicolumn{1}{c|}{212}                         & \multicolumn{1}{c|}{214}                         & \multicolumn{1}{c|}{0.94\%}                         & \multicolumn{1}{c|}{1800}                        & 33132         \\ \hline
\multicolumn{1}{|c|}{13}                         & \multicolumn{1}{c|}{17}                         & 227                     & \multicolumn{1}{c|}{225}                         & \multicolumn{1}{c|}{227}                         & \multicolumn{1}{c|}{0.89\%}                         & \multicolumn{1}{c|}{1800}                        & 2756          \\ \hline
\multicolumn{3}{|c|}{\textbf{Avg}}                                                                                           & \multicolumn{1}{c|}{138.4}                       & \multicolumn{1}{c|}{139.7}                       & \multicolumn{1}{c|}{0.98\%}                         & \multicolumn{1}{c|}{1041}                        & 538850        \\ \hline
\end{tabular}
\caption{Orbital branching performance summary.}
\label{table:orbital}
\end{table}

\begin{figure}[H]
\centering
    \includegraphics[scale=0.8]{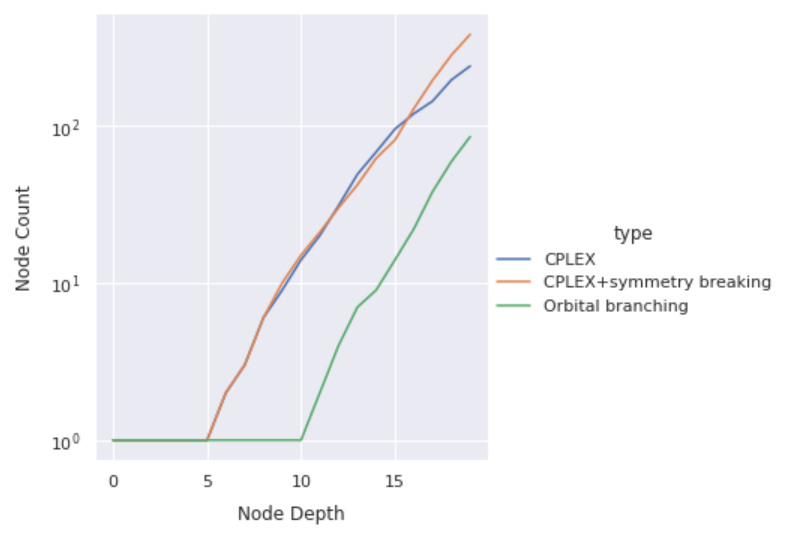}
    \caption{Branch-and-cut node count vs branch-and-cut node depth}
    \label{fig:orbitComparison}
\end{figure}

The results for Orbital Branching are shared in Table \ref{table:orbital}. Orbital Branching solves $9$ instances. Comparing Tables \ref{table:basiccp} and \ref{table:orbital}, we observe that while the Basic Formulation fails to find a solution to instances other than $m = Z(d)$, Orbital Branching solves all instances for $d = 7$ and $d = 8$. The difference can be attributed to how Orbital Branching prunes the initial nodes, which can be seen from Figure \ref{fig:orbitComparison}. Orbital Branching prunes all nodes but one until depth $10$ for the instance with $d=7$ and $m=8$. This enables it to work on a significantly smaller feasible region. Yet, when the problem size increases, its performance deteriorates. For $d=13$ and $m=17$, it fails to find a provably optimal solution, while CPLEX with aggressive symmetry breaking finds in fewer nodes. Some difference is due to nauty iterations running over $50$ seconds. Yet, for all methods, most instances are unsolved, and the problem of decreasing the upper bound continues.

In our next experiment, we explore the performance of Iterative Methods (Subsections \ref{iterative_section} and \ref{combined_section}). In the Iterative Method, we use the global callback mechanism of CPLEX. Instead of waiting for an iteration to end, we use the lower bound information from previous iterations to stop the current iteration if needed. We share the results for both  the Iterative Method and the Iterative Method with Orbital Branching in Table \ref{table:iterperf}. The Iterative Method optimally solves 14 instances, and the Iterative Method with Orbital Branching solves one more instance to optimality. These numbers are significantly higher than methods using Basic Formulation. The Iterative Method has slightly less run time, while the combined method solves the problem using significantly fewer nodes.

\begin{table}[]
\centering
\resizebox{\textwidth}{!}{\begin{tabular}{|ccc|ccccc|ccccc|}
\hline
\multicolumn{3}{|c|}{\textbf{Parameters}}                                                                                        & \multicolumn{5}{c|}{\textbf{Iterative}}                                                                                                                                                                                      & \multicolumn{5}{c|}{\textbf{Iterative + Orbital Branching}}                                                                                                                                                                                     \\ \hline
\multicolumn{1}{|c|}{\textit{\textbf{d}}}        & \multicolumn{1}{c|}{\textit{\textbf{m}}}        & \textbf{PreUB}              & \multicolumn{1}{c|}{\textbf{LB}}                 & \multicolumn{1}{c|}{\textbf{UB}}                 & \multicolumn{1}{c|}{\textbf{Gap}}                   & \multicolumn{1}{c|}{\textbf{Time}}               & \textbf{Node} & \multicolumn{1}{c|}{\textbf{LB}}                 & \multicolumn{1}{c|}{\textbf{UB}}                 & \multicolumn{1}{c|}{\textbf{Gap}}                   & \multicolumn{1}{c|}{\textbf{Time}}                & \textbf{Node}                   \\ \hline
\rowcolor[HTML]{C0C0C0} 
\multicolumn{1}{|c|}{\cellcolor[HTML]{C0C0C0}7}  & \multicolumn{1}{c|}{\cellcolor[HTML]{C0C0C0}8}  & 59                        & \multicolumn{1}{c|}{\cellcolor[HTML]{C0C0C0}58}  & \multicolumn{1}{c|}{\cellcolor[HTML]{C0C0C0}58}  & \multicolumn{1}{c|}{\cellcolor[HTML]{C0C0C0}0.00\%} & \multicolumn{1}{c|}{\cellcolor[HTML]{C0C0C0}0}   & 145           & \multicolumn{1}{c|}{\cellcolor[HTML]{C0C0C0}58}  & \multicolumn{1}{c|}{\cellcolor[HTML]{C0C0C0}58}  & \multicolumn{1}{c|}{\cellcolor[HTML]{C0C0C0}0.00\%} & \multicolumn{1}{c|}{\cellcolor[HTML]{C0C0C0}0}    & 57                              \\ \hline
\rowcolor[HTML]{C0C0C0} 
\multicolumn{1}{|c|}{\cellcolor[HTML]{C0C0C0}7}  & \multicolumn{1}{c|}{\cellcolor[HTML]{C0C0C0}9}  & 66                        & \multicolumn{1}{c|}{\cellcolor[HTML]{C0C0C0}66}  & \multicolumn{1}{c|}{\cellcolor[HTML]{C0C0C0}66}  & \multicolumn{1}{c|}{\cellcolor[HTML]{C0C0C0}0.00\%} & \multicolumn{1}{c|}{\cellcolor[HTML]{C0C0C0}0}   & 0             & \multicolumn{1}{c|}{\cellcolor[HTML]{C0C0C0}66}  & \multicolumn{1}{c|}{\cellcolor[HTML]{C0C0C0}66}  & \multicolumn{1}{c|}{\cellcolor[HTML]{C0C0C0}0.00\%} & \multicolumn{1}{c|}{\cellcolor[HTML]{C0C0C0}0}    & 0                               \\ \hline
\rowcolor[HTML]{C0C0C0} 
\multicolumn{1}{|c|}{\cellcolor[HTML]{C0C0C0}8}  & \multicolumn{1}{c|}{\cellcolor[HTML]{C0C0C0}9}  & 76                          & \multicolumn{1}{c|}{\cellcolor[HTML]{C0C0C0}74}  & \multicolumn{1}{c|}{\cellcolor[HTML]{C0C0C0}74}  & \multicolumn{1}{c|}{\cellcolor[HTML]{C0C0C0}0.00\%} & \multicolumn{1}{c|}{\cellcolor[HTML]{C0C0C0}3}   & 12880         & \multicolumn{1}{c|}{\cellcolor[HTML]{C0C0C0}74}  & \multicolumn{1}{c|}{\cellcolor[HTML]{C0C0C0}74}  & \multicolumn{1}{c|}{\cellcolor[HTML]{C0C0C0}0.00\%} & \multicolumn{1}{c|}{\cellcolor[HTML]{C0C0C0}4}    & 2794                            \\ \hline
\rowcolor[HTML]{C0C0C0} 
\multicolumn{1}{|c|}{\cellcolor[HTML]{C0C0C0}8}  & \multicolumn{1}{c|}{\cellcolor[HTML]{C0C0C0}10} & 84                          & \multicolumn{1}{c|}{\cellcolor[HTML]{C0C0C0}84}  & \multicolumn{1}{c|}{\cellcolor[HTML]{C0C0C0}84}  & \multicolumn{1}{c|}{\cellcolor[HTML]{C0C0C0}0.00\%} & \multicolumn{1}{c|}{\cellcolor[HTML]{C0C0C0}0}   & 0             & \multicolumn{1}{c|}{\cellcolor[HTML]{C0C0C0}84}  & \multicolumn{1}{c|}{\cellcolor[HTML]{C0C0C0}84}  & \multicolumn{1}{c|}{\cellcolor[HTML]{C0C0C0}0.00\%} & \multicolumn{1}{c|}{\cellcolor[HTML]{C0C0C0}0}    & 0                               \\ \hline
\rowcolor[HTML]{C0C0C0} 
\multicolumn{1}{|c|}{\cellcolor[HTML]{C0C0C0}9}  & \multicolumn{1}{c|}{\cellcolor[HTML]{C0C0C0}10} & 94                        & \multicolumn{1}{c|}{\cellcolor[HTML]{C0C0C0}92}  & \multicolumn{1}{c|}{\cellcolor[HTML]{C0C0C0}92}  & \multicolumn{1}{c|}{\cellcolor[HTML]{C0C0C0}0.00\%} & \multicolumn{1}{c|}{\cellcolor[HTML]{C0C0C0}2}   & 5149          & \multicolumn{1}{c|}{\cellcolor[HTML]{C0C0C0}92}  & \multicolumn{1}{c|}{\cellcolor[HTML]{C0C0C0}92}  & \multicolumn{1}{c|}{\cellcolor[HTML]{C0C0C0}0.00\%} & \multicolumn{1}{c|}{\cellcolor[HTML]{C0C0C0}11}   & 1990                            \\ \hline
\rowcolor[HTML]{C0C0C0} 
\multicolumn{1}{|c|}{\cellcolor[HTML]{C0C0C0}9}  & \multicolumn{1}{c|}{\cellcolor[HTML]{C0C0C0}11} & 103                       & \multicolumn{1}{c|}{\cellcolor[HTML]{C0C0C0}102} & \multicolumn{1}{c|}{\cellcolor[HTML]{C0C0C0}102} & \multicolumn{1}{c|}{\cellcolor[HTML]{C0C0C0}0.00\%} & \multicolumn{1}{c|}{\cellcolor[HTML]{C0C0C0}5}   & 12459         & \multicolumn{1}{c|}{\cellcolor[HTML]{C0C0C0}102} & \multicolumn{1}{c|}{\cellcolor[HTML]{C0C0C0}102} & \multicolumn{1}{c|}{\cellcolor[HTML]{C0C0C0}0.00\%} & \multicolumn{1}{c|}{\cellcolor[HTML]{C0C0C0}36}   & 4170                            \\ \hline
\rowcolor[HTML]{C0C0C0} 
\multicolumn{1}{|c|}{\cellcolor[HTML]{C0C0C0}9}  & \multicolumn{1}{c|}{\cellcolor[HTML]{C0C0C0}12} & 112                       & \multicolumn{1}{c|}{\cellcolor[HTML]{C0C0C0}112} & \multicolumn{1}{c|}{\cellcolor[HTML]{C0C0C0}112} & \multicolumn{1}{c|}{\cellcolor[HTML]{C0C0C0}0.00\%} & \multicolumn{1}{c|}{\cellcolor[HTML]{C0C0C0}1}   & 33            & \multicolumn{1}{c|}{\cellcolor[HTML]{C0C0C0}112} & \multicolumn{1}{c|}{\cellcolor[HTML]{C0C0C0}112} & \multicolumn{1}{c|}{\cellcolor[HTML]{C0C0C0}0.00\%} & \multicolumn{1}{c|}{\cellcolor[HTML]{C0C0C0}368}  & 229                             \\ \hline
\rowcolor[HTML]{C0C0C0} 
\multicolumn{1}{|c|}{\cellcolor[HTML]{C0C0C0}10} & \multicolumn{1}{c|}{\cellcolor[HTML]{C0C0C0}11} & 115                         & \multicolumn{1}{c|}{\cellcolor[HTML]{C0C0C0}112} & \multicolumn{1}{c|}{\cellcolor[HTML]{C0C0C0}112} & \multicolumn{1}{c|}{\cellcolor[HTML]{C0C0C0}0.00\%} & \multicolumn{1}{c|}{\cellcolor[HTML]{C0C0C0}100} & 208001        & \multicolumn{1}{c|}{\cellcolor[HTML]{C0C0C0}112} & \multicolumn{1}{c|}{\cellcolor[HTML]{C0C0C0}112} & \multicolumn{1}{c|}{\cellcolor[HTML]{C0C0C0}0.00\%} & \multicolumn{1}{c|}{\cellcolor[HTML]{C0C0C0}79}   & 10114                           \\ \hline
\rowcolor[HTML]{C0C0C0} 
\multicolumn{1}{|c|}{\cellcolor[HTML]{C0C0C0}10} & \multicolumn{1}{c|}{\cellcolor[HTML]{C0C0C0}12} & 125                         & \multicolumn{1}{c|}{\cellcolor[HTML]{C0C0C0}125} & \multicolumn{1}{c|}{\cellcolor[HTML]{C0C0C0}125} & \multicolumn{1}{c|}{\cellcolor[HTML]{C0C0C0}0.00\%} & \multicolumn{1}{c|}{\cellcolor[HTML]{C0C0C0}1}   & 0             & \multicolumn{1}{c|}{\cellcolor[HTML]{C0C0C0}125} & \multicolumn{1}{c|}{\cellcolor[HTML]{C0C0C0}125} & \multicolumn{1}{c|}{\cellcolor[HTML]{C0C0C0}0.00\%} & \multicolumn{1}{c|}{\cellcolor[HTML]{C0C0C0}2}    & 0                               \\ \hline
\rowcolor[HTML]{C0C0C0} 
\multicolumn{1}{|c|}{\cellcolor[HTML]{C0C0C0}11} & \multicolumn{1}{c|}{\cellcolor[HTML]{C0C0C0}12} & 137                       & \multicolumn{1}{c|}{\cellcolor[HTML]{C0C0C0}134} & \multicolumn{1}{c|}{\cellcolor[HTML]{C0C0C0}134} & \multicolumn{1}{c|}{\cellcolor[HTML]{C0C0C0}0.00\%} & \multicolumn{1}{c|}{\cellcolor[HTML]{C0C0C0}176} & 253162        & \multicolumn{1}{c|}{\cellcolor[HTML]{C0C0C0}134} & \multicolumn{1}{c|}{\cellcolor[HTML]{C0C0C0}134} & \multicolumn{1}{c|}{\cellcolor[HTML]{C0C0C0}0.00\%} & \multicolumn{1}{c|}{\cellcolor[HTML]{C0C0C0}146}  & 15106                           \\ \hline
\rowcolor[HTML]{C0C0C0} 
\multicolumn{1}{|c|}{\cellcolor[HTML]{C0C0C0}11} & \multicolumn{1}{c|}{\cellcolor[HTML]{C0C0C0}13} & 148                       & \multicolumn{1}{c|}{\cellcolor[HTML]{C0C0C0}146} & \multicolumn{1}{c|}{\cellcolor[HTML]{C0C0C0}146} & \multicolumn{1}{c|}{\cellcolor[HTML]{C0C0C0}0.00\%} & \multicolumn{1}{c|}{\cellcolor[HTML]{C0C0C0}206} & 362834        & \multicolumn{1}{c|}{\cellcolor[HTML]{C0C0C0}146} & \multicolumn{1}{c|}{\cellcolor[HTML]{C0C0C0}146} & \multicolumn{1}{c|}{\cellcolor[HTML]{C0C0C0}0.00\%} & \multicolumn{1}{c|}{\cellcolor[HTML]{C0C0C0}355}  & 371215                          \\ \hline
\multicolumn{1}{|c|}{11}                         & \multicolumn{1}{c|}{14}                         & 159                       & \multicolumn{1}{c|}{158}                         & \multicolumn{1}{c|}{159}                         & \multicolumn{1}{c|}{0.63\%}                         & \multicolumn{1}{c|}{1800}                        & 1087900       & \multicolumn{1}{c|}{158}                         & \multicolumn{1}{c|}{159}                         & \multicolumn{1}{c|}{0.63\%}                         & \multicolumn{1}{c|}{1800}                         & 706832                          \\ \hline
\rowcolor[HTML]{C0C0C0} 
\multicolumn{1}{|c|}{\cellcolor[HTML]{C0C0C0}11} & \multicolumn{1}{c|}{\cellcolor[HTML]{C0C0C0}15} & 170                       & \multicolumn{1}{c|}{\cellcolor[HTML]{C0C0C0}170} & \multicolumn{1}{c|}{\cellcolor[HTML]{C0C0C0}170} & \multicolumn{1}{c|}{\cellcolor[HTML]{C0C0C0}0.00\%} & \multicolumn{1}{c|}{\cellcolor[HTML]{C0C0C0}2}   & 38            & \multicolumn{1}{c|}{\cellcolor[HTML]{C0C0C0}170} & \multicolumn{1}{c|}{\cellcolor[HTML]{C0C0C0}170} & \multicolumn{1}{c|}{\cellcolor[HTML]{C0C0C0}0.00\%} & \multicolumn{1}{c|}{\cellcolor[HTML]{C0C0C0}1}    & 0                               \\ \hline
\multicolumn{1}{|c|}{12} & \multicolumn{1}{c|}{13} & 162 & \multicolumn{1}{c|}{158}                         & \multicolumn{1}{c|}{159}                         & \multicolumn{1}{c|}{0.63\%}                         & \multicolumn{1}{c|}{1800}                        & 2264281       & \multicolumn{1}{c|}{\cellcolor[HTML]{C0C0C0}158} & \multicolumn{1}{c|}{\cellcolor[HTML]{C0C0C0}158} & \multicolumn{1}{c|}{\cellcolor[HTML]{C0C0C0}0.00\%} & \multicolumn{1}{c|}{\cellcolor[HTML]{C0C0C0}1623} & \cellcolor[HTML]{C0C0C0}2020924 \\ \hline
\multicolumn{1}{|c|}{12}                         & \multicolumn{1}{c|}{14}                         & 174                         & \multicolumn{1}{c|}{171}                         & \multicolumn{1}{c|}{172}                         & \multicolumn{1}{c|}{0.58\%}                         & \multicolumn{1}{c|}{1800}                        & 1579790       & \multicolumn{1}{c|}{171}                         & \multicolumn{1}{c|}{172}                         & \multicolumn{1}{c|}{0.58\%}                         & \multicolumn{1}{c|}{1800}                         & 1101510                         \\ \hline
\rowcolor[HTML]{C0C0C0} 
\multicolumn{1}{|c|}{\cellcolor[HTML]{C0C0C0}12} & \multicolumn{1}{c|}{\cellcolor[HTML]{C0C0C0}15} & 186                         & \multicolumn{1}{c|}{\cellcolor[HTML]{C0C0C0}186} & \multicolumn{1}{c|}{\cellcolor[HTML]{C0C0C0}186} & \multicolumn{1}{c|}{\cellcolor[HTML]{C0C0C0}0.00\%} & \multicolumn{1}{c|}{\cellcolor[HTML]{C0C0C0}0}   & 0             & \multicolumn{1}{c|}{\cellcolor[HTML]{C0C0C0}186} & \multicolumn{1}{c|}{\cellcolor[HTML]{C0C0C0}186} & \multicolumn{1}{c|}{\cellcolor[HTML]{C0C0C0}0.00\%} & \multicolumn{1}{c|}{\cellcolor[HTML]{C0C0C0}0}    & 0                               \\ \hline
\multicolumn{1}{|c|}{13}                         & \multicolumn{1}{c|}{14}                         & 188                       & \multicolumn{1}{c|}{184}                         & \multicolumn{1}{c|}{185}                         & \multicolumn{1}{c|}{0.54\%}                         & \multicolumn{1}{c|}{1800}                        & 1495811       & \multicolumn{1}{c|}{184}                         & \multicolumn{1}{c|}{185}                         & \multicolumn{1}{c|}{0.54\%}                         & \multicolumn{1}{c|}{1800}                         & 2294463                         \\ \hline
\multicolumn{1}{|c|}{13}                         & \multicolumn{1}{c|}{15}                         & 201                       & \multicolumn{1}{c|}{198}                         & \multicolumn{1}{c|}{200}                         & \multicolumn{1}{c|}{1.01\%}                         & \multicolumn{1}{c|}{1800}                        & 1313790       & \multicolumn{1}{c|}{198}                         & \multicolumn{1}{c|}{200}                         & \multicolumn{1}{c|}{1.01\%}                         & \multicolumn{1}{c|}{1800}                         & 935875                          \\ \hline
\multicolumn{1}{|c|}{13}                         & \multicolumn{1}{c|}{16}                         & 214                       & \multicolumn{1}{c|}{212}                         & \multicolumn{1}{c|}{213}                         & \multicolumn{1}{c|}{0.47\%}                         & \multicolumn{1}{c|}{1800}                        & 690328        & \multicolumn{1}{c|}{212}                         & \multicolumn{1}{c|}{213}                         & \multicolumn{1}{c|}{0.47\%}                         & \multicolumn{1}{c|}{1800}                         & 577438                          \\ \hline
\rowcolor[HTML]{C0C0C0} 
\multicolumn{1}{|c|}{\cellcolor[HTML]{C0C0C0}13} & \multicolumn{1}{c|}{\cellcolor[HTML]{C0C0C0}17} & 227                       & \multicolumn{1}{c|}{\cellcolor[HTML]{C0C0C0}227} & \multicolumn{1}{c|}{\cellcolor[HTML]{C0C0C0}227} & \multicolumn{1}{c|}{\cellcolor[HTML]{C0C0C0}0.00\%} & \multicolumn{1}{c|}{\cellcolor[HTML]{C0C0C0}56}  & 5696          & \multicolumn{1}{c|}{\cellcolor[HTML]{C0C0C0}227} & \multicolumn{1}{c|}{\cellcolor[HTML]{C0C0C0}227} & \multicolumn{1}{c|}{\cellcolor[HTML]{C0C0C0}0.00\%} & \multicolumn{1}{c|}{\cellcolor[HTML]{C0C0C0}72}   & 543                             \\ \hline
\multicolumn{2}{|c|}{\textbf{Avg}}                                                                 & \multicolumn{1}{l|}{}       & \multicolumn{1}{c|}{138.45}                      & \multicolumn{1}{c|}{138.8}                       & \multicolumn{1}{c|}{0.25\%}                         & \multicolumn{1}{c|}{568}                         & 464615        & \multicolumn{1}{c|}{138.45}                      & \multicolumn{1}{c|}{138.75}                      & \multicolumn{1}{c|}{0.22\%}                         & \multicolumn{1}{c|}{585}                          & 402163                          \\ \hline
\end{tabular}}
\caption{Iterative Formulation performance summary.}
\label{table:iterperf}
\end{table}

\begin{table}[]
\centering
\begin{tabular}{|c|c|c|c|c|}
\hline
\textbf{Method}               & \textbf{Solved} & \textbf{Gap} & \textbf{Time} & \textbf{Node} \\ \hline
Basic Formulation                         & 7               & 1.12\%       & 1233          & 1099650       \\ \hline
Basic Formulation  + Symmetry Breaking     & 7               & 1.12\%       & 1207          & 729346        \\ \hline
Orbital Branching             & 9            & 0,98\%       & 1041          & 538850        \\ \hline
Iterative Method                  & 14              & 0,25\%       & 568           & 464615        \\ \hline
Iterative Method + Orbital Branching & 15              & 0.22\%       & 585           & 402163        \\ \hline
\end{tabular}
\caption{Performance summary of all methods.}
\label{table:summary}
\end{table}

The overall performance summary is shared in Table \ref{table:summary}. We can see an increase in instances solved as we use more elaborate approaches. The number of solved instances for the Basic Formulation is $7$, but all of them are for $m = Z(d)$ values; therefore, there is a significant increase in the performance with the new approaches. Especially, both  Iterative Methods are significantly better than the others.

\begin{figure}[htb]\label{fig:Bdd+t}
\begin{center}
\includegraphics[width=11cm]{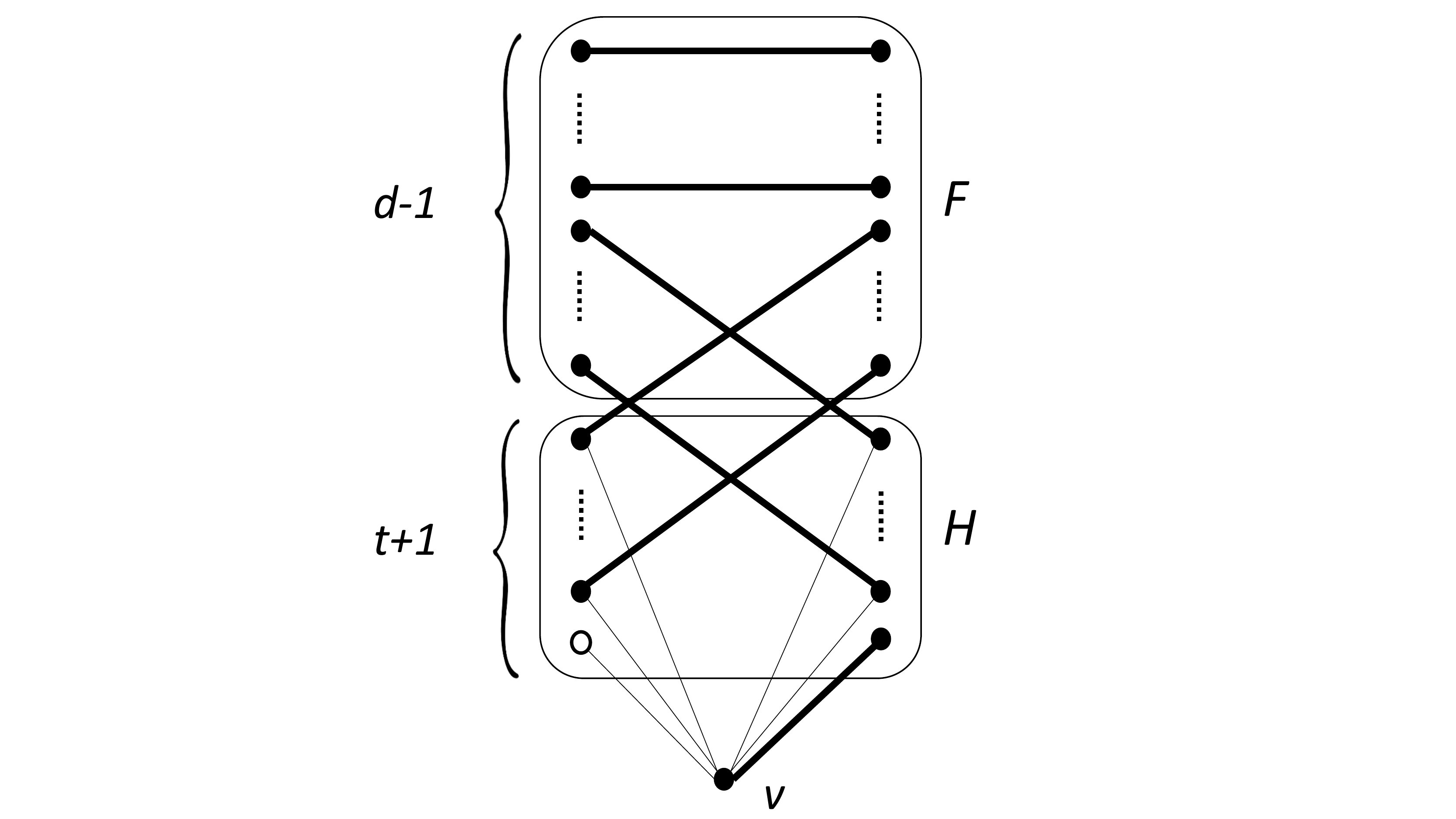}
\caption{The graph $B_{d,d+t}$ where the set $H$ induces an independent set and the set $F$ induces a complete bipartite graph from which $t$ distinct perfect matchings are removed. The deletion of the vertex $v$ leaves a bipartite graph that has complete links between the opposite parts of the sets $F$ and $H$. Bold edges show a maximum matching of size $d+t$.}
\end{center}
\end{figure}

It is important to note that Iterative Methods do not only solve some new instances that support Conjecture \ref{conj:comp}, but also provide us with graphs that follow a pattern. This trend allows us to suggest a general construction for graphs claimed to be extremal in Conjecture \ref{conj:comp} for all maximum degree $d$ and matching number $i$ with $7\leq d < i<Z(d)$. For $0\leq t< Z(d)-d$, let us describe the graph $B_{d,d+t}$, where $i = d + t$ as follows (see Figure \ref{fig:Bdd+t}): take a complete bipartite graph $K_{d+t,d+t}$ with $d+t$ vertices in each side. Consider $d-1$ vertices of each side (the set denoted $F$ in Figure \ref{fig:Bdd+t}) and remove $t$ disjoint perfect matchings between them. Note that this is possible because since $t< Z(d)-d$ and $Z(d)\le \lceil{5(d+1)/4}\rceil$ by Lemma \ref{theorem:z}, we have $t<\lceil{(d+1)/4}\rceil+1<d$. Then, consider the remaining $t+1$ vertices of each side (the set denoted $H$ in Figure \ref{fig:Bdd+t}), remove all edges between, introduce one additional vertex $v$, and make it adjacent to all the $t+1$ vertices of each side. With this construction, all the vertices except $v$ have degree $d$.
Let us remark that for $t=0$, the graph $B_{d,d+t}$ corresponds to the graph $A_d$ given in \cite{aey22} as an extremal graph for the case $d=m$ with $f_{\Delta}(d,d)=d^2+1$ edges. The graph $A_d$, or equivalently $B_{d,d}$, can be seen as a 5-cycle whose two adjacent vertices are replaced with independent sets of size $d-1$ each and edges between these sets and other vertices (or set) are replaced with complete links. From this perspective, our construction of the graph $B_{d,d+t}$ can be seen as a generalization of $A_d$ for matching numbers between $d+1$ and $Z(d)-1$.

\begin{proposition}\label{prop:newconstruction}
For all $i=d+t$ such that $d < i < Z(d)$, the graph $B_{d,d+t}$ has $di+i-d+1$ edges, 
$\Delta(B_{d,d+t})=d$ and $\nu({B_{d,d+t}})=d+t=i$; thus, it is an extremal graph if Conjecture \ref{conj:comp} holds.
\end{proposition}
\begin{proof}
First, we note that $t<Z(d)-d$ and thus the graph $B_{d,d+t}$ has one vertex of degree $2t+2$ (the additional vertex) and $2d+2t$ vertices each one of degree $d$. Since $t<\lceil{(d+1)/4}\rceil+1$, we have $2t+2\leq d$ for $d\geq 7$, it follows that $\Delta(B_{d,d+t})\leq d$. Moreover, we have $\nu({B_{d,d+t}})\leq d+t=i$ since it has $2d+2t+1$ vertices. More precisely a matching of size $d+t=i$ can be obtained as shown in Figure \ref{fig:Bdd+t} by bold edges; set a matching of size $d-t-2$ in $F$, match $v$ with one vertex in $H$, leave one vertex from the other side of $H$ unmatched, now since $t+1\leq d/2 <d-1$, all vertices of $H$ but these two can be matched with the remaining vertices of $F$ (recall that each side of $F$ is completely joined to the opposite side of $H$). Now, counting the number of edges using the degrees yields $d^2+dt+t+1$ edges which is, for $i=d+t$, equal to $di+i-d+1$ as suggested in Conjecture \ref{conj:comp}. 
\end{proof}

\begin{table}[h]
\centering
\begin{tabular}{|c|c|c|c|c|c|}
\hline
\multicolumn{1}{|l|}{$\boldsymbol{d}$} & \multicolumn{1}{l|}{$\boldsymbol{m}$} & \multicolumn{1}{l|}{\textbf{edge count}} & \multicolumn{1}{l|}{\textbf{\#$\boldsymbol{d}$-star}} & \multicolumn{1}{l|}{\#$\boldsymbol{Z(d)}$} & \multicolumn{1}{l|}{\textbf{\#other}} \\ \hline
7                                & 15                              & 108                                      & 6                                    & 1                                    & 0                                     \\ \hline
7                                & 16                              & 116                                      & 0                                    & 1                                    & 1                                     \\ \hline
7                                & 17                              & 124                                      & 0                                    & 1                                    & 1                                     \\ \hline
7                                & 18                              & 132                                      & 0                                    & 2                                    & 0                                     \\ \hline
7                                & 19                              & 139                                      & 1                                    & 2                                    & 0                                     \\ \hline
7                                & 20                              & 146                                      & 2                                    & 2                                    & 0                                     \\ \hline
7                                & 21                              & 153                                      & 3                                    & 2                                    & 0                                     \\ \hline
8                                & 17                              & 140                                      & 7                                    & 1                                    & 0                                     \\ \hline
8                                & 18                              & 148                                      & 8                                    & 1                                    & 0                                     \\ \hline
8                                & 19                              & 158                                      & 0                                    & 1                                    & 1                                     \\ \hline
8                                & 20                              & 168                                      & 0                                    & 2                                    & 0                                     \\ \hline
8                                & 21                              & 176                                      & 1                                    & 2                                    & 0                                     \\ \hline
8                                & 22                              & 184                                      & 2                                    & 2                                    & 0                                     \\ \hline
8                                & 23                              & 192                                      & 3                                    & 2                                    & 0                                     \\ \hline
8                                & 24                              & 200                                      & 4                                    & 2                                    & 0                                     \\ \hline
9                                & 19                              & 175                                      & 7                                    & 1                                    & 0                                     \\ \hline
9                                & 20                              & 184                                      & 8                                    & 1                                    & 0                                     \\ \hline
9                                & 21                              & 194                                      & 0                                    & 1                                    & 1                                     \\ \hline
9                                & 22                              & 204                                      & 0                                    & 1                                    & 1                                     \\ \hline
9                                & 23                              & 214                                      & 0                                    & 1                                    & 1                                     \\ \hline
9                                & 24                              & 224                                      & 0                                    & 2                                    & 0                                     \\ \hline
9                                & 25                              & 233                                      & 1                                    & 2                                    & 0                                     \\ \hline
9                                & 26                              & 242                                      & 2                                    & 2                                    & 0                                     \\ \hline
9                                & 27                              & 251                                      & 3                                    & 2                                    & 0                                     \\ \hline
10                               & 21                              & 215                                      & 9                                    & 1                                    & 0                                     \\ \hline
10                               & 22                              & 226                                      & 0                                    & 1                                    & 1                                     \\ \hline
10                               & 23                              & 237                                      & 0                                    & 1                                    & 1                                     \\ \hline
10                               & 24                              & 250                                      & 0                                    & 2                                    & 0                                     \\ \hline
10                               & 25                              & 260                                      & 1                                    & 2                                    & 0                                     \\ \hline
10                               & 26                              & 270                                      & 2                                    & 2                                    & 0                                     \\ \hline
10                               & 27                              & 280                                      & 3                                    & 2                                    & 0                                     \\ \hline
10                               & 28                              & 290                                      & 4                                    & 2                                    & 0                                     \\ \hline
10                               & 29                              & 300                                      & 5                                    & 2                                    & 0                                     \\ \hline
10                               & 30                              & 310                                      & 6                                    & 2                                    & 0                                     \\ \hline
\end{tabular}
\caption{Knapsack Formulation results.}
\label{table:count}
\end{table}

Finally, we use the above mentioned results to construct an extremal graphs for all $d$ values $7, 8, 9, 10$ with the Knapsack Formulation. For these $d$ values, there are at most $4$ components in the Knapsack Formulation; therefore, even for $d = 10$ and $m = 10000$, it takes 0.2 seconds to obtain an optimal solution. The number of each extremal component in an optimal solution and resulting edge counts are shared in Table \ref{table:count}. In Table \ref{table:count}, edge counts are shared for $2d < m \le 3d$. The column ``\textit{\#star}" denotes the number of $d$-stars, ``\textit{\#Z(d)}" denotes the number of graphs with matching number $Z(d)$, and ``\textit{\#other}" denotes the number of other components. In Table \ref{table:count8}, more extensive results are shared for $d=8$ where ``\textit{comp\_i}" denotes the number of components with matching number $i$ in an extremal graph. When $m$ increases by $10$, the only difference is an additional component with matching number $Z(8)=10$.  Tables \ref{table:count} and \ref{table:count8} show that Proposition \ref{prop:number}, which is conditional to Conjecture \ref{conj:comp}, holds for $d=7,8,9$ and 10. Proposition \ref{prop:number} describes a pattern on the construction of an edge-extremal graph; it suggests that there is an extremal graph having as many extremal components as possible with matching number $Z(d)$, completed with at most one extremal component with matching number between $d$ and $Z(d)$, and $d$-stars; all of our results are compatible with this foresight and support therefore both Conjecture \ref{conj:comp} and \ref{conj:lemma}. It follows that the edge numbers of extremal graphs for all $m$ values and for $d \in \{7,8,9,10\}$ are equal to the formula given in Conjecture \ref{conj:lemma}.

\begin{table}[H]
\centering
\begin{tabular}{|c|c|c|c|c|c|c|}
\hline
$\boldsymbol{d}$ & $\boldsymbol{m}$ & \textbf{edge count} & \textbf{$\boldsymbol{d}$-star} & \textbf{comp\_8} & \textbf{comp\_9} & \textbf{comp\_10} \\ \hline
8          & 15         & 124                 & 5                & 0                & 0                & 1                 \\ \hline
8          & 16         & 132                 & 6                & 0                & 0                & 1                 \\ \hline
8          & 17         & 140                 & 7                & 0                & 0                & 1                 \\ \hline
8          & 18         & 149                 & 0                & 1                & 0                & 1                 \\ \hline
8          & 19         & 158                 & 0                & 0                & 1                & 1                 \\ \hline
8          & 20         & 168                 & 0                & 0                & 0                & 2                 \\ \hline
8          & 21         & 176                 & 1                & 0                & 0                & 2                 \\ \hline
8          & 22         & 184                 & 2                & 0                & 0                & 2                 \\ \hline
8          & 23         & 192                 & 3                & 0                & 0                & 2                 \\ \hline
8          & 24         & 200                 & 4                & 0                & 0                & 2                 \\ \hline
8          & 25         & 208                 & 5                & 0                & 0                & 2                 \\ \hline
8          & 26         & 216                 & 6                & 0                & 0                & 2                 \\ \hline
8          & 27         & 224                 & 7                & 0                & 0                & 2                 \\ \hline
8          & 28         & 233                 & 0                & 1                & 0                & 2                 \\ \hline
8          & 29         & 242                 & 0                & 0                & 1                & 2                 \\ \hline
8          & 30         & 252                 & 0                & 0                & 0                & 3                 \\ \hline
8          & 31         & 260                 & 1                & 0                & 0                & 3                 \\ \hline
8          & 32         & 268                 & 2                & 0                & 0                & 3                 \\ \hline
8          & 33         & 276                 & 3                & 0                & 0                & 3                 \\ \hline
8          & 34         & 284                 & 4                & 0                & 0                & 3                 \\ \hline
8          & 35         & 292                 & 5                & 0                & 0                & 3                 \\ \hline
8          & 36         & 300                 & 6                & 0                & 0                & 3                 \\ \hline
8          & 37         & 308                 & 7                & 0                & 0                & 3                 \\ \hline
8          & 38         & 317                 & 0                & 1                & 0                & 3                 \\ \hline
8          & 39         & 326                 & 0                & 0                & 1                & 3                 \\ \hline
8          & 40         & 336                 & 0                & 0                & 0                & 4                 \\ \hline
8          & 41         & 344                 & 1                & 0                & 0                & 4                 \\ \hline
8          & 42         & 352                 & 2                & 0                & 0                & 4                 \\ \hline
8          & 43         & 360                 & 3                & 0                & 0                & 4                 \\ \hline
8          & 44         & 368                 & 4                & 0                & 0                & 4                 \\ \hline
8          & 45         & 376                 & 5                & 0                & 0                & 4                 \\ \hline
8          & 46         & 384                 & 6                & 0                & 0                & 4                 \\ \hline
8          & 47         & 392                 & 7                & 0                & 0                & 4                 \\ \hline
\end{tabular}
\caption{Knapsack Formulation results for $d=8$.}
\label{table:count8}
\end{table}

\section{Conclusion}
In this paper, we studied the open cases for the problem of finding the maximum number of edges in a triangle-free graph with maximum degree at most $d$ and matching number at most $m$. We suggested several integer programming methods, which is, to the best of our knowledge, a new approach for this extremal problem.

Since our Basic Formulation is highly symmetric, most of our efforts have been about breaking the symmetry. We proposed two different approaches, first exploiting symmetry for branching decisions using Orbital Branching and then an Iterative Method exploiting the closeness of the precomputed upper bound and the optimal solution value. As we used a combination of the Iterative Method and Orbital Branching, we observed an increase in the number of instances optimally solved and an overall decrease in solution time. Our approach extends the known cases from $d\leq 6$ in \cite{aey22} to $d\leq 10$ (both for all $m>d$). We could also identify some extremal components for $d=11, 12$ and $13$. The only missing extremal component for $d=11$ is for $m=14$; if this extremal component is found then the Knapsack Formulation can compute all extremal graphs for $d=11$ and any $m>d$. 

It is important to note that although integer programming approaches will be limited in finding extremal components for higher $d$ and $m$ values, they give us more evidence on Conjectures \ref{conj:comp}, \ref{conj:lemma} and Proposition \ref{prop:number}. Indeed, all our findings support the claim that there is an extremal triangle-free graph with as many factor-critical extremal components with matching number $Z(d)$ as possible , at most one factor-critical extremal component with matching number less than $Z(d)$, and some $d$-stars. Thus, our results provide additional motivation to search for a formal (structural) proof of Conjectures \ref{conj:comp} and \ref{conj:lemma}, which will most probably require more powerful tools. 

Finally, we pose the following question. We note that our integer programming formulations exploit the fact that extremal components are factor-critical by using the implied vertex number in their formulations. However, they do not explicitly force the constructed graphs to be factor-critical. Clearly, a factor-critical graph with matching number $m$ has $2m+1$ vertices, but a graph with $2m+1$ vertices is not necessarily factor-critical. However, it turns out that all extremal components resulting from our integer programming formulations are factor-critical. This observation suggests that it might be interesting to investigate the following question: is it true that a triangle-free extremal graph $G$ with matching number $m$ and maximum degree $d$ such that $d < m < Z(d)$ and having $2m+1$ vertices is factor-critical?

\section*{Acknowledgments}
This work has been supported by the Scientific and Technological Research Council of Turkey (T\"UB\.ITAK) under the grant number 122M452.
%We thank to the anonymous referees for their valuable suggestions which improved the presentation of our paper. 

\FloatBarrier
\printbibliography %Prints bibliography

@article{chem1,
author = {Hansen, P. and M. Aouchiche and G. Caporossi and A. Hertz and C. Sellal},
journal = {Journal of Chemistry and Chemical Engineering Systems},
pages = {22-30},
title = {Mixed Integer Programming and Extremal Chemical Graphs},
volume = {3},
year = {2018},
}

@book{chem2,
author = {Ito, R. and N. A. Azam and C. Wang and A. Shurbevski and H. Nagamochi and T. Akutsu},
journal = {Advances in Computer Vision and Computational Biology},
pages = {641-655},
title = {A novel method for the inverse QSAR/QSPR to monocyclic chemical compounds based on artificial neural networks and integer programming},
publisher = {Springer},
address = {Cham},
year = {2021},
}

@article{furini,
author = {Furini, Fabio and Ljubi\'c, Ivana and Segundo, Pablo San},
journal = {Manuscript},
number = {},
pages = {},
title = {A New Bilevel Optimization Approach for Computing {Ramsey} Numbers},
volume = {},
year = {2021},
URL = {https://optimization-online.org/?p=17323}
}

@techreport{aey22,
author = {Ahanjideh, Milad and T{\i}naz Ekim and Mehmet Akif Y{\i}ld{\i}z},
institution = {arXiv},
title = {Maximum size of a triangle-free graph with bounded maximum degree and matching number},
type = {preprint},
year = {2022},
archivePrefix = {arXiv},
eprint = {2207.02271},
}

@article{bk09,
author = {Balachandran, Niranjan and Niraj Khare},
journal = {Discrete Mathematics},
number = {12},
pages = {4176-4180},
title = {Graphs with restricted valency and matching number},
volume = {309},
year = {2009},
}

@article{b22,
author = {Blair, J.R.S. and Heggernes, P. and Lima, P.T. and Lokshtanov, D.},
journal = {Algorithmica}, 
pages = {3587-3602},
title = {On the maximum number of edges in chordal graphs of bounded degree and matching number},
volume = {84},
year = {2022},
}

@article{ch76,
author = {Chv{\'a}tal, V. and D. Hanson},
journal = {Journal of Combinatorial Theory, Series B},
number = {2},
pages = {128-138},
title = {Degrees and matchings},
volume = {20},
year = {1976},
}

@article{deh17,
author = {Dibek, Cemil and T{\i}naz Ekim and Pinar Heggernes},
journal = {Discrete Mathematics},
number = {5},
pages = {927-934},
title = {Maximum number of edges in claw-free graphs whose maximum degree and matching number are bounded},
volume = {340},
year = {2017},
}

@article{er60,
author = {Erdös, Paul and Richard Rado},
journal = {Journal of the London Mathematical Society},
number = {1},
pages = {85-90},
title = {Intersection theorems for systems of sets},
volume = {1},
year = {1960},
}

@article{e64,
author = {Turan, P.},
journal = {Középisk. Mat. és Fiz. Lapok},
pages = {436-452},
title = {On an extremal problem in graph theory},
volume = {48},
year = {1941},
}

@mastersthesis{m15,
author = {M{\aa}land, Erik Kvam},
school = {The University of Bergen},
title = {Maximum number of edges in graph classes under degree and matching constraints},
year = {2015},
}

@article{m03,
author = {Margot, Fran\c{c}ois},
journal = {Mathematical Programming},
number = {1},
pages = {3-21},
title = {Exploiting orbits in symmetric {ILP}},
volume = {98},
year = {2003},
}

@article{mp14,
author = {McKay, Brendan D. and Adolfo Piperno},
journal = {Journal of Symbolic Computation},
pages = {94-112},
title = {Practical graph isomorphism, II},
volume = {60},
year = {2014},
}

@article{j05,
author = {Johnson, David S.},
journal = {ACM Transactions on Algorithms (TALG)},
number = {1},
pages = {160-176},
title = {The {NP}-completeness column},
volume = {1},
year = {2005},
}

@article{pr19,
author = {Pfetsch, Marc E. and Thomas Rehn},
journal = {Mathematical Programming Computation},
number = {1},
pages = {37-93},
title = {A computational comparison of symmetry handling methods for mixed integer programs},
volume = {11},
year = {2019},
}

@article{o11,
author = {Ostrowski, J. and Linderoth, J. and Rossi, F. and Smriglio, S.},
journal = {Mathematical Programming},
number = {1},
pages = {147-178},
title = {Orbital branching},
volume = {126},
year = {2011},
}

\end{document}